\documentclass[12pt]{article} 
\usepackage[margin=1in]{geometry}
\usepackage{setspace}
\usepackage{graphicx}	
\usepackage{subcaption}
\usepackage{float} 
\usepackage[toc]{appendix}
\usepackage[draft]{todonotes} 
\usepackage{enumerate}	
\usepackage{paralist}	
\usepackage{amsmath, amsthm, amssymb,mathtools}
\usepackage{thm-restate}

\usepackage[pdftex,hypertexnames=false,linktocpage=true,hidelinks]{hyperref}
\usepackage{cleveref}

\newtheorem{thm}{Theorem}[section]

\newtheorem{prop}[thm]{Proposition}

\newtheorem{lemma}[thm]{Lemma}

\theoremstyle{definition}

\theoremstyle{definition}

\usepackage{tikz}
\usepackage{xcolor}

\usetikzlibrary{decorations.pathreplacing}

\definecolor{amber}{rgb}{1.0, 0.75, 0.0}

\definecolor{forest}{rgb}{0.0, 0.5, 0.0}

\definecolor{cadmium}{rgb}{0.93, 0.53, 0.18}

\definecolor{byzantine}{rgb}{0.74, 0.2, 0.64}

\definecolor{brilliantrose}{rgb}{1.0, 0.33, 0.64}

\definecolor{caribbeangreen}{rgb}{0.0, 0.8, 0.6}

\definecolor{electriccyan}{rgb}{0.0, 1.0, 1.0}

\definecolor{periwinkle}{rgb}{0.8, 0.8, 1.0}

\definecolor{steelblue}{rgb}{0.27, 0.51, 0.71}

\title{Rainbow Tur\'an numbers for short brooms}

\author{John Byrne\footnote{Department of Mathematical Sciences, University of Delaware, Newark, DE, USA. E-mail: \texttt{jpbyrne@udel.edu}. Research partially supported by NSF DSM-2245556.} \and E.G.K.M. Gamlath\footnote{College of Health and Natural Sciences, Franklin Pierce University, Rindge, NH.  E-mail: \texttt{egkmgamlath@gmail.com}.} \and Anastasia Halfpap\footnote{Department of Mathematics, Iowa State University, Ames, IA, USA.  E-mail: \texttt{ahalfpap@iastate.edu}. Research supported by NSF DSM-2152490.} \and Sydney Miyasaki\footnote{Department of Mathematics, Iowa State University, Ames, IA, USA.  E-mail: \texttt{miyasaki@iastate.edu}. Research supported by NSF DMS-2152490} \and Alex Parker\footnote{Department of Mathematics, Iowa State University, Ames, IA, USA.  E-mail: \texttt{abparker@iastate.edu}.}}

\begin{document}
\maketitle

\begin{abstract}
A graph $G$ is \textit{rainbow}-$F$-\textit{free} if it admits a proper edge-coloring without a rainbow copy of $F$. The \textit{rainbow Tur\'an number} of $F$, denoted $\mathrm{ex^*}(n,F)$, is the maximum number of edges in a rainbow-$F$-free graph on $n$ vertices. We determine bounds on the rainbow Tur\'an numbers of stars with a single edge subdivided twice; we call such a tree with $t$ total edges a $t$-edge \textit{broom} with length-$3$ handle, denoted by $B_{t,3}$. 
We improve the best known upper bounds on $\mathrm{ex^*}(n,B_{t,3})$ in all cases where $t \neq 2^s - 2$. Moreover, in the case where $t$ is odd and in some cases when $t \equiv 0 \mod 4$, we provide constructions asymptotically achieving these upper bounds. Our results also demonstrate a dependence of $\mathrm{ex^*}(n,B_{t,3})$ on divisibility properties of $t$.

\end{abstract}

\section{Introduction}

Given a fixed graph $F$, we say that $G$ \textit{contains} $F$ if some subgraph of $G$ (not necessarily induced) is isomorphic to $F$. We call a subgraph of $G$ isomorphic to $F$ an $F$-\textit{copy} in $G$. If $G$ does not contain $F$ as a subgraph, we say that $G$ is $F$-\textit{free}. A foundational question in extremal combinatorics is to determine the maximum number of edges in an $n$-vertex $F$-free graph. We denote this maximum by $\mathrm{ex}(n,F)$, the \textit{Tur\'an number} of $F$. Since the problem of determining of $\mathrm{ex}(n,F)$ was initiated by Mantel and Tur\'an, this function has been very broadly studied, and has inspired a number of variations and generalizations. 

An \textit{edge-coloring} of a graph $G$ is a function $c: E(G) \rightarrow \mathbb{N}$, where we call $c(e)$ the \textit{color} of edge $e$. We say that an edge-coloring $c$ is \textit{proper} if for any incident edges $e,f$, we have $c(e) \neq c(f)$. An edge-coloring $c$ is \textit{rainbow} if $c$ is injective, i.e., for any edges $e,f$, we have $c(e) \neq c(f)$. Given fixed graphs $G,F$ and an edge-coloring $c$ of $G$, we say that $G$ is \textit{rainbow-F-free} if no $F$-copy in $G$ is rainbow under $c$. The \textit{rainbow Tur\'an number} of $F$, $\mathrm{ex}^*(n,F)$, is the maximum number of edges in an $n$-vertex graph $G$ which admits a proper edge-coloring under which $G$ is rainbow-$F$-free. Study of $\mathrm{ex}^*(n,F)$ was initiated by Keevash, Mubayi, Sudakov, and Verstra\"ete~\cite{KMSV} in 2007, and has received substantial attention since its introduction. 

It is clear that for any $F$, we have $\mathrm{ex}(n,F) \leq \mathrm{ex}^*(n,F)$, since any proper edge-coloring of an $n$-vertex, $F$-free graph is rainbow-$F$-free. In fact, $\mathrm{ex}(n,F)$ and $\mathrm{ex}^*(n,F)$ can only differ by a sub-quadratic term.

\begin{thm}[Keevash-Mubayi-Sudakov-Verstra\"ete~\cite{KMSV}]
For any graph $F$,
\[\mathrm{ex}(n,F) \leq \mathrm{ex}^*(n,F) \leq \mathrm{ex}(n,F) + o(n^2).\]
\end{thm}

Since $\mathrm{ex}(n,F) = \Theta(n^2)$ when $F$ is non-bipartite~\cite{ErSi,ErSt}, it follows that $\mathrm{ex}(n,F)$ and $\mathrm{ex}^*(n,F)$ are often asymptotically equal. However, for bipartite $F$, we often have that $\mathrm{ex}^*(n,F)$ differs from $\mathrm{ex}(n,F)$ by a constant multiplicative factor. It seems difficult in general to exactly determine $\mathrm{ex}^*(n,F)$ when $F$ is bipartite, even in cases where $\mathrm{ex}(n,F)$ is well-understood. 
As an example of these contrasts, let $P_{\ell}$ denote the $\ell$-edge path. We have $\mathrm{ex}(n,P_{\ell}) = \frac{\ell - 1}{2}n + O(1)$ for all $\ell$ (see~\cite{ErdosGallai}), but there is no general formula known for $\mathrm{ex}^*(n,P_{\ell})$. By a construction of Johnston and Rombach~\cite{JoRo}, we know that $\mathrm{ex^*}(n,P_{\ell}) \geq \frac{\ell}{2}n +O(1)$ for all $\ell \geq 3$, and thus $\mathrm{ex}(n,P_{\ell})$ is not asymptotically equal to $\mathrm{ex}^*(n,P_{\ell})$ in general. When $\ell \in \{3,4,5\}$, we have $\mathrm{ex}^*(n,P_{\ell}) = \frac{\ell}{2}n + O(1)$ (see \cite{halfpaprainbowp5,JoPaSa}), but the correct leading coefficient for $\mathrm{ex}^*(n,P_{\ell})$ is not known when $\ell \geq 6$. 

Paths are far from the only trees whose rainbow Tur\'an numbers are not well-understood. Another example is $B_{t,\ell}$, the \textit{broom} with $t$ edges and handle length $\ell$; this is the graph obtained from a \textit{handle} $P_\ell$ by adding $t-\ell$ pendant edges to an end vertex; these additional pendant edges are called \textit{bristles}. Note that $B_{t,t} = P_t$ and $B_{t+1,t} = P_{t+1}$. Brooms are a natural family of trees to consider in their own right, but they can also be viewed as a ``midpoint'' between stars (whose rainbow Tur\'an numbers are trivial to compute) and paths (whose rainbow Tur\'an numbers are very difficult to understand). Johnston and Rombach~\cite{JoRo} previously studied rainbow Tur\'an numbers for brooms (among other trees), and exactly characterized $\mathrm{ex^*}(n,B_{t,2})$ for all $t$.

\begin{thm}[Johnston-Rombach~\cite{JoRo}]
\[\mathrm{ex^*}(n,B_{t,2}) = \begin{cases} \frac{t}{2}n + O(1) & \text{ for $t$ odd} \\ \frac{t^2}{2(t+1)}n + O(1) & \text{ for $t$ even}\end{cases}\]
\end{thm}

Johnston and Rombach also determined $\mathrm{ex^*}(n,B_{t,3})$ exactly for certain values of $t$.

\begin{thm}[Johnston-Rombach~\cite{JoRo}]\label{special even t}
If $t = 2^s - 2$ for some $s \geq 3$, then
\[\mathrm{ex^*}(n,B_{t,3}) = \frac{t+1}{2}n + O(1).\]
\end{thm}

No other exact results are known on rainbow Tur\'an numbers for brooms. We note that other authors have considered rainbow Tur\'an numbers for trees more generally; in particular, Bednar and Bushaw \cite{bednar2022rainbow} give some general bounds, including an upper bound for caterpillars (a class of tree which includes all brooms). It seems to become more difficult to compute $\mathrm{ex^*}(n,B_{t,\ell})$ exactly as $\ell$ grows; the larger $\ell$ is, the more ``path-like'' $B_{t,\ell}$ becomes. For a more general upper bound when $\ell$ is not too large in comparison to $t$, Johnston and Rombach~\cite{JoRo} showed that
\[\mathrm{ex^*}(n,B_{t,\ell}) \leq \frac{t + \ell - 2}{2}n + O(1)\]
as long as $3\ell - 4 \leq t$. When $\ell$ is very small, this general bound is sometimes known to be tight: it gives the correct value when $\ell = 2$ and $t$ is odd, and when $\ell = 3$ and $t = 2^s - 2$ for some $s \geq 3$. However, as we shall see, there are also many cases where the above bound can be improved.

We build upon these results by studying the rainbow Tur\'an numbers of brooms with handles of length 3. We are able to precisely determine $\mathrm{ex^*}(n,B_{t,3})$ when $t$ is odd.

\begin{restatable}{thm}{oddtupperbound}\label{odd t upper bound}
Let $t \geq 3$ be an odd integer. Then $$\mathrm{ex}^*(n, B_{t,3}) = \frac{t}{2}n + O(1)$$.

\end{restatable}

The behavior of $\mathrm{ex^*}(n,B_{t,3})$ is somewhat more complicated when $t$ is even. Theorem~\ref{special even t} determines $\mathrm{ex^*}(n,B_{t,3})$ when $t$ is an even number of the form $2^s - 2$. If $t$ is not of the form $2^s - 2$, we show that $\mathrm{ex^*}(n,B_{t,3})$ is strictly smaller than $\frac{t+1}{2}n$.

\begin{restatable}{thm}{eventupperbound}\label{even t upper bound}
Let $t$ be an even integer such that $t\ne 2^s-2$ for any $s$. Then $$\frac{t-1}{2}n+O(1)\le\mathrm{ex}^*(n,B_{t,3})\le\left(\frac{t+1}{2}-\frac{1}{t+2}\right)n.$$

\end{restatable}

When $t\equiv 0\pmod 4$, we prove a better upper bound.

\begin{restatable}{thm}{multipleoffourupperbound}\label{multiple of 4 upper bound}
Let $t\equiv 0\pmod 4$. Then $$\mathrm{ex}^*(n,B_{t,3})\le\frac{t}{2}n+O(1).$$

\end{restatable}

When $t$ is a power of two, we prove a matching lower bound. 

\begin{restatable}{thm}{poweroftwoextremalnumber}\label{power of 2 extremal number}
Let $t=2^s$ for some $s$. Then $$\mathrm{ex}^*(n,B_{t,3})=\frac{t}{2}n+O(1).$$

\end{restatable}

Finally, we show that the lower bound in Theorem~\ref{power of 2 extremal number} also holds when $t$ is one less than a power of $3$. Note that $3^s - 1 \equiv 0 \pmod{4}$ if and only if $s$ is even, so Theorem~\ref{multiple of 4 upper bound} and Theorem~\ref{power of 3 lower bound} combine to determine $\mathrm{ex}^*(n,B_{t,3})$ in half of the cases where $t = 3^s - 1$.

\begin{restatable}{thm}{powerofthreelowerbound}\label{power of 3 lower bound}
Let $t=3^s-1$ for some $s\ge 2$. Then $$\mathrm{ex}^*(n,B_{t,3}) \ge \frac{t}{2}n + O(1).$$

\end{restatable}

The remainder of the paper is organized as follows. In Section~\ref{odd broom section}, we prove Theorem~\ref{odd t upper bound}, as well as some lemmas which are useful throughout the paper. In Section~\ref{even section}, we prove the theorems concerning even values of $t$. In Section~\ref{conclusion}, we give concluding remarks.

\subsection{Notation and definitions}

Throughout, we use standard notation. We write $[n]$ for the set of integers $\{1,\ldots,n\}$. Given a graph $G$ and $v \in V(G)$, we denote by $d(v)$ the degree of $v$. We denote by $\delta(G)$ and $\overline{d}(G)$ the minimum and average degrees, respectively, of $G$. The \textit{neighborhood} of $v$, denoted $N(v)$, is the set $\{u \in V(G): uv \in E(G)\}$, while the \textit{closed neighborhood} of $v$, denoted $N[v]$, is $N(v) \cup \{v\}$. Given an edge-coloring $c:E(G)\to\mathbb N$, we write $\mathrm{im}(c)$ for $\{c(e):e\in E(G)\}$.

Constructions in this paper will be described in terms of several common graphs and subgraphs. We denote by $K_k$ the clique on $k$ vertices, and by $K_{s,t}$ the complete bipartite graph with class sizes $s$ and $t$. In a graph $G$, a \textit{matching} $M$ is a subset of $E(G)$ such that no vertex of $G$ is contained in more than one element of $M$. We say that $M$ is a \textit{perfect} matching if $M$ is a matching and every vertex of $G$ is contained in some edge of $M$.

\section{Odd brooms}\label{odd broom section}

\subsection{Lower bound}

\begin{thm}\label{broom lower bound}
If $t$ is odd, then we have $\mathrm{ex}^*(n,B_{t,3})\ge\frac{t}{2}n+O(1)$.
\end{thm}
\begin{proof}
    First we show that $K_{t+1}$ admits a proper edge-coloring with no rainbow $B_{t,3}$. Since $t+1$ is even, it is possible to partition the edges of $K_{t+1}$ into perfect matchings (see e.g.~\cite{Kirkman}). Let each perfect matching correspond to a unique color in $\{1,\ldots,t\}$. Suppose that a rainbow $B_{t,3}$ appears, with handle $v_1,v_2,v_3,v_4$ and bristles $v_5,\ldots,v_{k+1}$. Assume that $c(v_1v_2)=1$. Since the broom has an edge between $v_4$ and every vertex of $\{v_3\}\cup\{v_5,\ldots,v_{k+1}\}$, the edge in color $1$ which is incident to $v_4$ must be $v_1v_4$ or $v_2v_4$; but in either case, the coloring is not proper.

    Now for $n\in\mathbb{N}$, taking $\left\lfloor\frac{n}{k+1}\right\rfloor$ disjoint copies of $K_{k+1}$ along with some isolated vertices and coloring each clique as described above gives
    $$\mathrm{ex}^*(n,B_{k,3})\ge\left\lfloor\frac{n}{k+1}\right\rfloor\frac{(k+1)k}{2}=\frac{k}{2}n+O(1).$$
\end{proof}

\subsection{Upper bound}

By Theorem~\ref{broom lower bound}, we have $\mathrm{ex}^*(n,B_{t,3}) \geq \frac{t}{2} n + O(1)$ for any odd $t$. In this section we provide a matching upper bound to prove Theorem~\ref{odd t upper bound}, which we restate here for convenience.

\oddtupperbound*

When $t = 3$, we have $B_{3,3} = P_3$, and we know that $\mathrm{ex}^*(n,P_3) = \frac{3}{2}n + O(1)$ \cite{JoPaSa}. Thus, we focus on the cases where $t > 3$. 

Before proving Theorem~\ref{odd t upper bound}, we shall require some set-up and a few lemmas. Our broad strategy is to argue by contradiction. If Theorem~\ref{odd t upper bound} is untrue, we will be able to find a graph $G^*$ with several properties, and obtain a contradiction by analyzing $G^*$. The existence of this $G^*$ is guaranteed by the following folklore lemma. 

\begin{lemma}\label{good subgraph}
Let $G$ be a graph with $\overline{d}(G) \geq d$. Then $G$ contains a subgraph $G^*$ such that $\delta(G^*) > \frac{d}{2}$ and $\overline{d}(G^*) \geq d$.
\end{lemma}

Fix $t \geq 5$ an odd integer, and suppose there exists some $n$-vertex graph $G$ with $e(G) > \frac{t}{2}n$ which admits a rainbow-$B_{t,3}$-free proper edge-coloring. By Lemma~\ref{good subgraph}, $G$ contains a subgraph $G^*$ with $\delta(G^*) > \frac{t}{2}$ and $\overline{d}(G^*) > t$. We may moreover assume that every component of $G^*$ has average degree strictly larger than $t$ (since if not, components with low average degree may be deleted). For the remainder of this section, we will work with $G^*$.

The following lemma applies to $G^*$, and will help to determine the structure of $G^*$ around a fixed vertex.

\begin{lemma}[Halfpap~\cite{halfpaprainbowp5}]\label{P3 lemma}
Suppose $G$ is a graph such that $\delta(G) \geq 3$ and every component of $G$ has average degree at least $5$. Let $v \in V(G)$. Then $v$ is the endpoint of a rainbow $P_3$-copy.
\end{lemma}

We now prove some lemmas concerning the structure of $G^*$. We begin with the following simple structural lemma which will apply to the components of $G^*$ and also be useful in Section~\ref{even section}. 

\begin{lemma}\label{extremal space restriction}
Let $t \geq 5$ and suppose $G$ is a connected graph with $\overline{d}(G) > t$ and $\delta(G) \geq 3$. If $G$ admits a rainbow-$B_{t,3}$-free proper edge-coloring, then $G$ is a subgraph of $K_{t+2}$.
\end{lemma}

\begin{proof}

Let $v$ be a vertex of degree (at least) $t+1$ in $G$. By Lemma~\ref{P3 lemma}, there exists a $P_3$-copy $P$ with $v$ as an endpoint, say $P = vxyz$, with $c(vx) = 1$, $c(xy) = 2$, and $c(yz) = 3$. Now, $v$ has at least $t-2$ neighbors which are not on $P$, say $v_1, v_2, \dots v_{t-2}$. To avoid a rainbow $B_{t,3}$-copy immediately, we have (without loss of generality) $c(vv_i) = i+1$ for each $i \in [t-2]$. This also implies that to avoid a rainbow $B_{t,3}$-copy, $v$ must have exactly $t-2$ neighbors which are not on $P$. So $vy, vz \in E(G)$; without loss of generality, $c(vy) = t$ and $c(vz) = t+1$. We illustrate the configuration in the case $t = 8$ in Figure~\ref{8,3 subgraph}.

\begin{figure}[h]

\begin{center}

\begin{tikzpicture}

\draw[thick, red] (0,0) -- (2,0) node[pos = 0.5, above]{1};
\draw[thick, cadmium] (2,0) -- (4,0) node[pos = 0.5, above]{2};
\draw[thick, amber] (4,0) -- (6,0) node[pos = 0.5, above]{3};
\draw[thick, cadmium] (0,0) -- (-1.25,1.5) node[pos = 0.9, right]{2};
\draw[thick, amber] (0,0) -- (-1.75,1);
\draw[amber] (-1.4,1) node{3};
\draw[thick, forest] (0,0) -- (-2.25,0.5);
\draw[forest] (-1.9,0.62)node{4};
\draw[thick, blue] (0,0) -- (-2.25,-0.5);
\draw[blue] (-1.9,-0.62)node{5};
\draw[thick, brilliantrose] (0,0) -- (-1.25,-1.5) node[pos = 0.9, right]{7};
\draw[thick, byzantine] (0,0) -- (-1.75,-1);
\draw[byzantine] (-1.4,-1) node{6};

\draw[thick, cyan] (0,0) to[bend left = 50] (4,0);
\draw[cyan] (2,1.15) node{8};
\draw[thick, gray] (0,0) to[bend right = 40] (6,0);
\draw[gray] (3,-1.4) node{9};

\filldraw (0,0) circle (0.05 cm) node[above]{$v$};
\filldraw (2,0) circle (0.05 cm)node[above]{$x$};
\filldraw (4,0) circle (0.05 cm)node[above]{$y$};
\filldraw (6,0) circle (0.05 cm)node[above]{$z$};
\filldraw (-1.75,1) circle (0.05 cm)node[left]{$v_2$};
\filldraw (-1.75,-1) circle (0.05 cm)node[left]{$v_5$};
\filldraw (-2.25,0.5) circle (0.05 cm)node[left]{$v_3$};
\filldraw (-2.25,-0.5) circle (0.05 cm)node[left]{$v_4$};
\filldraw (-1.25,1.5) circle (0.05 cm)node[left]{$v_1$};
\filldraw (-1.25,-1.5) circle (0.05 cm)node[left]{$v_6$};

\end{tikzpicture}

\caption{The structure of $N[v]$}\label{8,3 subgraph}

\end{center}

\end{figure}
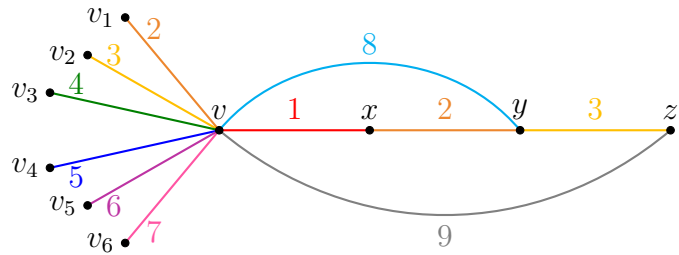

We now show that for every $u \in N[v]$, we have $N[u] \subseteq N[v]$.  
Note that by the above analysis, if there exists a rainbow $P_3$-copy $P$ with $v$ as an endpoint, then $v$ must be adjacent to every vertex on $P$ (excepting itself). Observe that if some $u \in N(v)$ is adjacent to $w \not \in N[v]$, then  $vuw$ is a rainbow $P_2$-copy ending at $v$ such that $v$ is not adjacent to $w$, so there must be no neighbor $s$ of $w$ such that $vuws$ is rainbow. In particular, this means that  $|N(w) \setminus \{v,u\}| \leq 1$ (since any neighbor of $w$ not on $vuw$ must be adjacent to $w$ via an edge of color $c(vu)$). This implies $d(w) \leq 2$, a contradiction. We conclude that there does not exist any path $vuw$ with $w \not \in N(v)$, so $N[u] \subseteq N[v]$.
\end{proof}

Lemma~\ref{extremal space restriction} implies that each component of $G^*$ contains at most $t+2$ vertices. Since $\overline{d}(G^*)\ge\overline{d}(G)$, this implies that either $\mathrm{ex^*}(n,B_{t,3}) \leq \frac{t}{2}n$ or an extremal construction for $\mathrm{ex^*}(n,B_{t,3})$ is provided by taking disjoint copies of some (possibly not proper) subgraph of $K_{t+2}$ with an appropriate edge-coloring.

We also note that Lemma~\ref{extremal space restriction} applies for any value of $t \geq 5$, not only odd $t$, so we shall also apply it in Section~\ref{even section} to bound $\mathrm{ex^*}(n,B_{t,3})$ for even $t$.

Next, we work to bound the size of common neighborhoods in $G^*$.
Note that if $N(v) \cap N(u)$ is large, then there exist many $C_4$-copies containing both $v$ and $u$. The following observation will be useful in understanding the colorings of such $C_4$-copies.

\begin{prop}\label{c4 colorings}

Fix $t\geq 3$ and suppose $G$ is a properly edge-colored, rainbow-$B_{t,3}$-free graph. Suppose $vxyzv$ is a $C_4$-copy in $G$ and $|N(v) \setminus \{x,y,z\}| \geq t-2$. Then either $vxyzv$ is bichromatic, or $vxyzv$ is rainbow and $v$ is incident to edges of colors $c(xy)$ and $c(yz)$.

\end{prop}
\begin{proof}
Suppose for a contradiction that $vxyzv$ is a $C_4$-copy whose coloring is not among the described types. First, suppose $vxyzv$ is colored with precisely $3$-edge colors. Without loss of generality, $c(vx) = c(yz) =1$, $c(xy) = 2$, and $c(zv) = 3$. At most one vertex in $N(v) \setminus \{x,y,z\}$ is adjacent to $v$ via an edge whose color is in $\{1,2,3\}$; thus, at least $t-3$ vertices in $N(v)\setminus \{x,y,z\}$ are adjacent to $v$ via edges whose colors are not from $\{1,2,3\}$. Thus, $N[v]$ 
contains a rainbow $B_{t,3}$ with handle $xyzv$, a contradiction.

Next, suppose $vxyzv$ is a rainbow $C_4$, say with $c(vx) = 1$, $c(xy) = 2$, $c(yz) = 3$, $c(zv) = 4$, and $v$ is not incident to an edge colored $3$. As above, at most one vertex in $N(v) \setminus \{x,y,z\}$ is adjacent to $v$ via an edge whose color is in $\{1,2,3,4\}$; thus, at least $t-3$ vertices in $N(v)\setminus \{x,y,z\}$ are adjacent to $v$ via edges whose colors are not from $\{1,2,3,4\}$. Thus, $N[v]$ contains a rainbow $B_{t,3}$ with handle $xyzv$, a contradiction.
\end{proof}

We note that, like Lemma~\ref{extremal space restriction}, Proposition~\ref{c4 colorings} applies for any value of $t\geq 3$, not just odd $t$. We will also use Proposition~\ref{c4 colorings} in Section~\ref{even section}.

We now consider the intersections of neighborhoods in $G^*$.

\begin{lemma}\label{derangements}
For any $u,v \in V(G^*)$, $|N(u) \cap N(v)| < t$. 
\end{lemma}

\begin{proof}
Suppose for a contradiction that there are some $u,v \in V(G^*)$ with $|N(u) \cap N(v)| \geq t$. Let $w_1, \dots, w_t$ be $t$ vertices in $N(u) \cap N(v)$. We define
\[C_u = \{c(uw_i): 1 \leq i \leq t\}\]
and 
\[C_v = \{c(vw_i): 1 \leq i \leq t\}.\]
We claim first that $C_u = C_v$. If not, consider $a \in C_v \setminus C_u$. Without loss of generality, $a = c(vw_1)$. Consider the $C_4$-copy $u w_1 v w_2 u$. Say $b = c(u w_1)$, and $c = c(u w_2)$. By Proposition~\ref{c4 colorings}, $c(vw_2) \neq b$, so $c(vw_2)$ must be a new color, say $d$. Now, using Proposition~\ref{c4 colorings} again, $u$ must be incident to an edge of color $a$. However, we know that this edge is not of the form $uw_i$ for any $i \in [t]$, since $a \not\in C_u$. Since the component of $u$ has at most $t+2$ vertices by Lemma~\ref{extremal space restriction}, the only possible neighbor of $u$ which is not in $\{w_i\}_{i=1}^t$ is $v$. Since $c(vw_1) = a$, if $c(uv) = a$ then the coloring is not proper, which is a contradiction.

Thus, we may assume $C_u = C_v$. Without loss of generality, $c(uw_i) = i$. Let $\sigma$ be the derangement of $[t]$ such that $c(vw_i) = \sigma(i)$. Observe, $u w_i v w_{\sigma(i)} u$ contains two edges of color $\sigma(i)$, so by Proposition~\ref{c4 colorings} must be bichromatic. Thus, we must have $c(vw_{\sigma(i)}) = i$, that is, $\sigma(\sigma(i)) = i$. We conclude that to avoid a rainbow $B_{t,3}$-copy, $\sigma$ must be a product of disjoint transpositions; however, since $t$ is odd, this is not possible.
\end{proof}

Using Lemma~\ref{derangements}, we derive a somewhat stronger statement about the interaction between large closed neighborhoods $N[u]$ with  vertices adjacent to $u$.

\begin{lemma}\label{closed neighborhood}
Suppose $u,v \in V(G^*)$ and $d(u) = t+1$. Then $|N[u] \cap N(v)| < t$.
\end{lemma}

\begin{proof}
Assume for the sake of a contradiction that $|N[u] \cap N(v)| \geq t$. Firstly, note that since components in $G^*$ have size $t+2$, if $d(u) = t+1$, then either $|N(u) \cap N(v)| = 0$ (that is, $u$ and $v$ are in different components) or $N(v) \subset N[u]$.

By Lemma~\ref{derangements}, we cannot have $|N(u) \cap N(v)| = t$. Thus, we may assume $|N(u) \cap N(v)| = t-1$; say $N(u) \cap N(v) = \{w_1, \dots, w_{t-1}\}$. Since $d(u) = t+1$, $u$ has exactly two neighbors which are not in $N(v)$, namely $v$ itself and another vertex, say $x$. Similarly to the proof of Lemma~\ref{derangements}, let 
\[C_u = \{c(uw_i): 1 \leq i \leq t-1\}\]
and
\[C_v = \{c(vw_i): 1 \leq i \leq t-1\}.\]
We now have two cases, depending upon the interaction of $C_u$ and $C_v$. 

\begin{enumerate}

\item $C_u \neq C_v$.

Without loss of generality, $c(uw_i) = i$ for each $i \in [t-1]$, and $c(vw_1) = t$. 
Observe, Proposition~\ref{c4 colorings} implies that if $c(vw_i) \not\in C_u$ for some $i \in [t]$, then every $C_4$-copy $u w_i v w_j u$ is rainbow, and moreover that $u$ is incident to an edge colored $c(vw_i)$. In particular, if $c(v w_i)\notin C_u$, then we must have $c(ux) = c(vw_i)$ (since the only neighbors of $u$ not in $\{w_i\}_{i=1}^{t-1}$ are $v$ and $x$, and we cannot have $c(uv) = c(vw_i)$). This implies that $C_v\setminus\{t\} \subset C_u$. 

Now, observe that if $c(v w_j) = 1$ for any $j \in \{2, \dots, t-1\}$, then $u w_1 v w_j u$ is a $3$-colored $C_4$-copy containing $v$, which does not exist by Proposition~\ref{c4 colorings}. Thus, $1 \not\in C_v$. We conclude that $C_v \setminus \{t\} = C_u \setminus \{1\}$, and in particular, 
\[C_v' := \{c(v w_i): 2 \leq i \leq t-1\}\]
is equal to $\{2, \dots, t-1\}$. Now, as in the proof of Lemma~\ref{derangements}, let $\sigma$ be the derangement of $\{2, \dots, t-1\}$ such that $c(v w_i) = \sigma(i)$. To avoid a $3$-colored $C_4$-copy, $\sigma$ must be a product of disjoint transpositions, for otherwise we can find $i,j\in\{2,\ldots,t-1\}$ such that $c(w_iv)=j$, and $c(vw_j)\ne i$; then $w_jvw_iu$ forms the handle of a $B_{t,3}$. However, since $t$ is odd and $\sigma$ is a permutation of $t-2$ objects, this is not possible.

\item $C_u = C_v$.

Without loss of generality, $C_u = C_v = [t-1]$, $c(uv) = t$, and $c(ux) = t+1$. We examine the possible neighbors of $x$. Recall that $x \not\in N(v)$. 

Observe that $x$ is not adjacent to any $w_i$. Indeed, if $x w_i$ is an edge and we have $c(x w_i) \in C_v$, then $w_i x u v$ is the handle of a rainbow $B_{t,3}$-copy, since $v$ has $t-2$ neighbors in $\{w_{j}: j \neq i\}$, of which only one is adjacent to $v$ via an edge of color $c(xw_i)$. On the other hand, if $c(x w_i) \not \in C_v$, then $u x w_i v u$ is either a trichromatic $C_4$-copy (if $c(x w_i) = t$) or a rainbow $C_4$-copy such that $u$ is not incident to an edge of color $c(x w_i)$; both outcomes are precluded by Proposition~\ref{c4 colorings}.

Since components in $G^*$ have size $t+2$ by Lemma~\ref{extremal space restriction}, we must have 
\[N(x) \subseteq \{u,v\} \cup \{w_i\}_{i=1}^{t-1}.\] Thus, the only possible neighbors of $x$ are $u$ and $v$, so $d(x) \leq 2$,  a contradiction, as $\delta(G^*) > \frac{t}{2} > 2$. \qedhere
\end{enumerate}
\end{proof}

With Lemmas~\ref{extremal space restriction} and \ref{closed neighborhood} in hand, we can now quickly establish Theorem~\ref{odd t upper bound}.

\begin{proof}[Proof of Theorem~\ref{odd t upper bound}]

When $t = 3$, we have $B_{3,3} = P_3$, and $\mathrm{ex^*}(n,P_3) = \frac{3}{2}n + O(1)$ (see~\cite{JoPaSa}). Thus, we suppose $t \geq 5$. As previously in this section, we suppose for a contradiction that $G^*$ is an $n$-vertex graph with $e(G^*) > \frac{t}{2}n$. Recall, we may assume that $\delta(G^*) > \frac{t}{2}$, that $\overline{d}(G^*) > t$, and that every component of $G^*$ has average degree greater than $t$. By Lemma~\ref{extremal space restriction}, each component of $G^*$ is a subgraph of $K_{t+2}$. Thus, consider a component of $G^*$, say $C$. By the above properties, $C$ contains a vertex $u$ with $d(u) = t+1$, and $N[u] = V(C)$.

Now by Lemma~\ref{closed neighborhood}, for any $v \in V(C)$ with $v \neq u$, we have $|N[u] \cap N(v)| < t$, and by the properties of $G^*$,
\[|N[u] \cap N(v)| = |N(v)| = d(v).\]
We conclude that $u$ is the only vertex in $C$ of degree $t+1$, and all other vertices in $C$ in fact have degree strictly smaller than $t$. Thus, $\overline{d}(C) < t$, a contradiction.
\end{proof}

\section{Even brooms}\label{even section}

\subsection{General upper bound}\label{even broom upper sec}

Existing results give the asymptotics of $\mathrm{ex}^*(n,B_{t,3})$ for even $t$ only when $t=2^s-2$ for some $s$, according to Theorem~\ref{special even t}. In this subsection we prove Theorem~\ref{even t upper bound}, showing that the asymptotic behavior of $\mathrm{ex}^*(n,B_{t,3})$ in fact depends on whether $t+2$ is a power of 2.

\eventupperbound*

The lower bound can be obtained by arbitrarily coloring copies of $K_t$, and in the rest of this subsection we prove the upper bound. In~\cite{JoRo} an edge-coloring of $K_{2^s}$ was constructed which shows that if $t=2^s-2$ then $\mathrm{ex}^*(n,B_{t,3})\ge\frac{t+1}{2}n+O(1)$. We will show that if such a coloring of $K_{t+2}$ exists, then $t+2$ must be a power of 2.

\begin{prop} \label{powers of 2}
    Suppose $K_{t+2}$ is properly edge-colored with no rainbow $B_{t,3}$. Then $t+2=2^s$ for some $s$.
\end{prop}
\begin{proof}
    Let $c$ be such an edge-coloring of $K_{t+2}$. First we show that $c$ must be \textit{optimal}, i.e. every color is a perfect matching of $K_{t+2}$. If $c$ is not optimal, the without loss of generality color 1 is not a perfect matching. By the upper bound of Theorem~\ref{odd t upper bound}, $t$ must be even, so there are two vertices $u,v$ not incident to color 1. Let $x,y$ be vertices such that $c(xy)=1$. Without loss of generality, $c(yu)=2$ and $c(uv)=3$. Now $v$ has $t-2$ neighbors outside $\{x,y,u\}$ which are not in color 1 or 3, so at least $t-3$ of them are not in color 1, 2, or 3. Thus, $xyuv$ is the handle of a rainbow $B_{t,3}$, a contradiction.
    
    Thus, assume $c$ is optimal. By Proposition~\ref{c4 colorings}, we have the following `4-cycle property': for any $x,y,z,w\in V(K_{t+2})$, $c(xy)=c(zw)$ implies $c(yz)=c(xw)$. Let $s$ be the largest integer $i$ such that we can find a set $U$ of $2^i$ vertices in $K_{t+2}$ for which only $2^i-1$ colors appear in $U$, i.e. the restriction of $c$ to $U$ partitions the edges in $U$ into perfect matchings. If $U=V(K_{t+2})$, we are done. Otherwise, let $y\in V(K_{t+2})-U$. Fix some $x\in U$, and note that $c(xy)$ does not appear in $U$. Since every vertex of $U$ is incident to $c(xy)$ and the coloring is proper, we find a set $U'$ such that $|U|=|U'|$, $U\cap U'=\emptyset$, and there is a matching between $U$ and $U'$ in color $c(xy)$, corresponding to a bijection $f:U\to U'$. By the 4-cycle property, we have $c(f(u)f(v))=c(uv)$ for every $u,v\in U$. Thus, $\{c(uv):u,v\in U\}=\{c(u'v'):u',v'\in U'\}$. Now let $u\in U$ and $u'\in U'$. For each $v\in U-u$, there is a unique $v'\in U'-u'$ such that $c(uv)=c(u'v')$, and this implies $c(vv')=c(uu')$. Hence, there is a matching between $U$ and $U'$ in color $c(uu')$. This shows that the edges between $U$ and $U'$ decompose into perfect matchings. In total there are $2^n-1+2^n=2^{n+1}-1$ colors appearing in $U\cup U'$, contradicting the maximality of $U$.
\end{proof}

In fact, this proof also gives the uniqueness of the optimal coloring constructed by Johnston and Rombach. Take $s$ to be the largest $i$ such that there is a set $U$ of $2^i$ vertices such that the restriction of $c$ to $U$ can be obtained by some identification $\phi:U\to\mathbb F_2^i$ so that $c(uv)=\phi(u)-\phi(v)$. Identify $c(xy)$ with the $(i+1)^{st}$ unit vector $e_{i+1}$ in $\mathbb F_2^{i+1}$, and extend the identification to $U\cup U'$ so that $\phi(f(u))=\phi(u)+e_{i+1}$. For any $uv\in U$, the 4-cycle property gives $c(f(u)f(v))=c(uv)=\phi(u)-\phi(v)=(\phi(u)+e_{i+1})-(\phi(v)+e_{i+1})=\phi(f(u))-\phi(f(v))$. For each of the $2^s-1$ colors besides $e_{i+1}$ which appear in the matchings between $U$ and $U'$, identify the color $c(xu')$ with $\phi(u')-\phi(x)$. Now let $v\in U$, $v'\in U'$. Let $u'$ be such that $\phi(u')=\phi(v')+\phi(x)-\phi(v)$. By the above, this implies $c(xv)=c(u'v')$, so by the 4-cycle property, $c(vv')=c(xu')=\phi(x)-\phi(u')=\phi(x)-\phi(v')-\phi(x)+\phi(v)=\phi(v')-\phi(v)$. This contradicts the maximality of $U$.

\begin{proof}[Proof of Theorem~\ref{even t upper bound}]
Suppose $t\ne 2^s-2$ and $G$ is an edge-colored graph with no rainbow $B_{t,3}$. By Lemma~\ref{extremal space restriction}, each component of $G$ with average degree greater than $t$ is a subgraph of $K_{t+2}$. By Proposition~\ref{powers of 2}, each such component of $G$ is then a proper subgraph of $K_{t+2}$, hence has average degree at most $t+1-\frac{2}{t+2}$. So, the average degree of $G$ is at most $t+1-\frac{2}{t+2}$.
\end{proof}

\subsection{Upper bound for $t \equiv 0 \pmod 4$}

When $t$ is even, the formula for $\mathrm{ex^*}(n,B_{t,3})$ will vary with the value of $t$. Indeed, Theorem~\ref{even t upper bound} tells us that $\mathrm{ex^*}(n,B_{t,3}) = \frac{t+1}{2}n + O(1)$ if and only if $t$ is of the form $2^s - 2$. If $t$ is not of the form $2^s - 2$, it is not clear where $\mathrm{ex^*}(n,B_{t,3})$ should fall within the range given by Theorem~\ref{even t upper bound}. Towards answering this question, we can improve the upper bound of Theorem~\ref{even t upper bound} in the case when $t \equiv 0 \pmod 4$. The goal of this subsection is to prove Theorem~\ref{multiple of 4 upper bound}, restated here for convenience.

\multipleoffourupperbound*

To prove Theorem~\ref{multiple of 4 upper bound}, we need the following structural lemma. We state it here and will prove it shortly. 
\begin{lemma}\label{structure}
Let $G$ be a rainbow $B_{t,3}$-free graph on at most $t+2$ vertices, where $t \equiv 0 \pmod{4}$. Then $G$ has at most two vertices of degree $t+1$. Additionally, if $G$ does have two such vertices, then every other vertex has degree at most $t-1$.
\end{lemma}
With Lemma~\ref{structure}, Theorem~\ref{multiple of 4 upper bound} follows easily.

\begin{proof}[Proof of Theorem~\ref{multiple of 4 upper bound}]
Consider a rainbow $B_{t,3}$-free graph $G$ on $n$ vertices and with $m = \mathrm{ex}^*(n,B_{t,3})$ edges. Note that if $t = 4$, then $B_{t,3} = P_4$, which does satisfy Theorem~\ref{multiple of 4 upper bound} (see~\cite{JoPaSa}). As such, we may assume $t > 4$. Suppose, for the sake of contradiction, that $m > \frac{t}{2}n$. Then $\overline{d}(G) > t$ and so there is a component $C$ with $\overline{d}(C) > t$. Using Lemma~\ref{good subgraph}, we may assume $\delta(G) \geq 3$ after passing to a subgraph, if necessary. By Lemma~\ref{extremal space restriction}, $C$ is a graph on at most $t+2$ vertices. Since $\overline{d}(C) > t$, $C$ has at least one vertex of degree $t+1$. By Lemma~\ref{structure}, it has at most two such vertices.  We now consider the cases.
In the first case, suppose $C$ has exactly one vertex of degree $t+1$. Note that $t+1$ is odd and that $C$ must have an even number of odd-degree vertices. Hence, some vertex in $C$ must have a degree which is at most $t-1$. Thus,
$$
\overline{d}(C) \leq \frac{t+1 + t - 1 + t^2}{t+2} = \frac{t^2 + 2t}{t+2} = t.
$$
In the second case, suppose $C$ has exactly two vertices of degree $t+1$. By Lemma~\ref{structure}, all other vertices in $C$ have degree $t - 1$, at most. Thus,
$$
\overline{d}(C) \leq \frac{2(t+1) + t(t-1)}{t+2} = \frac{t^2 + t + 2}{t+2} = t - 1 + \frac{4}{t+2}  \leq t.$$
In either case, we have the contradiction $t < \overline{d}(C) \leq t$ . Hence, $m \leq \frac{t}{2} n$.
\end{proof}
Now we need to prove Lemma~\ref{structure}. To do so, we will need some notation and terminology.

Let $t \in \mathbb{N}$ and let $H$ be an arbitrary graph on at most $t+2$ vertices and let $c$ be a proper edge-coloring of $H$. For any two vertices $u$ and $v$ in $H$, we let the sets $E_u, E_v \subset E(H)$ be the sets of edges in $H$ incident to $u$ and $v$, respectively. Define the subgraph $H_{u,v}$ of $H$ by $H_{u,v} = (V(H), E_u \cup E_v)$.  We define the function
\[\sigma_{u,v}: \{c(uw): w \in N_H(u) \cap N_H(v)\} \to \{c(vw): w \in N_H(u) \cap N_H(v)\}\] 
by $\sigma_{u,v}(c(uw)) = c(vw).$  We also say that $c$ is a \textit{good} coloring of $H$ if $c$ uses at most $t+1$ colors and every four-cycle in $H$ is either bichromatic or rainbow under $c$.

 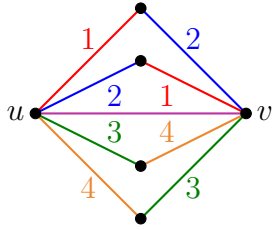
\begin{figure}[h]
 \begin{center}
 \begin{tikzpicture}[scale=1.4]
 \draw[thick, red] (0,0) -- (1,1) node[pos=0.5, above]{1};
 \draw[thick, blue] (0,0) -- (1,0.5) node[pos=0.75, below]{2};
 \draw[thick, forest] (0,0) -- (1,-0.5) node[pos=0.75, above]{3};
 \draw[thick, cadmium] (0,0) -- (1,-1) node[pos=0.5, below]{4};
 \draw[thick, blue] (2,0) -- (1,1) node[pos=0.5, above]{2};
 \draw[thick, red] (2,0) -- (1,0.5) node[pos=0.75, below]{1};
 \draw[thick,cadmium] (2,0) -- (1,-0.5) node[pos=0.75, above]{4};
 \draw[thick, forest] (2,0) -- (1,-1) node[pos=0.5, below]{3};
 \draw[thick, byzantine] (0,0) -- (2,0);
 \filldraw (0,0) circle (0.05 cm) node[left]{$u$};
 \filldraw (2,0) circle (0.05 cm) node[right]{$v$};
 \filldraw (1,1) circle (0.05 cm) ;
 \filldraw (1,0.5) circle (0.05 cm);
 \filldraw (1,-0.5) circle (0.05 cm);
 \filldraw (1,-1) circle (0.05 cm);
 \end{tikzpicture}
 \caption{$H_{u,v}$ where $t=4$, with good coloring $c$. Here $\sigma_{u,v} = (12)(34)$}
 \end{center}
 \end{figure}

Good colorings interact nicely with the functions $\sigma_{u,v}$, which follows the same argument as in Lemma~\ref{derangements}. 
\begin{lemma}\label{high degree derangements}
Let $H$ be a graph on at most $t + 2$ vertices with a good coloring c. If $u$ and $v$ are vertices of degree $t+1$, then $\sigma_{u,v}$ is both a derangement of $\mathrm{im}(c) \setminus \{c(uv)\}$ and a product of disjoint transpositions.
\end{lemma}

Now we take a fixed graph $G$ on at most $t+2$ vertices and a proper edge-coloring $c$ under which $G$ is rainbow-$B_{t,3}$-free. The connection between the structure of our fixed graph $G$ and these definitions in expressed in the following lemma, which follows from Proposition~\ref{c4 colorings}.

\begin{lemma}\label{good restriction}
If $u$ is a vertex of degree $t+1$ in $G$ and $c$ is a rainbow-$B_{t,3}$-free coloring of $G$, then $c$ restricts to a good coloring of $G_{u,v}$ for any vertex $v$.
\end{lemma}
The crux of our argument lies in the relationships between these derangements:
\begin{lemma}\label{high degree in good}
Let $H$ be a graph on at most $t+2$ vertices with a good coloring $c$. If $H$ has three vertices of degree $t+1$, then $t \equiv 2 \pmod{4}$. 
\end{lemma}
\begin{proof}
Let $x,y$ and $z$ be three vertices in $H$ of degree $t+1$. Without loss of generality, we may assume that the range of $c$ is $[t+1]$, that $c(xz) = 1$, $c(yz) = 2$, and that $c(xy)= t+1$. Let $v$ be the unique vertex satisfying $c(xv) = 2$. We note that, by the goodness of $H$, we have $c(yv) = 1$, and by the goodness of $H$ and the colors on the cycle $xzvyx$, we have $c(zv) = t+1$. Hence, $(12)$ is a transposition in the disjoint cycle decomposition of $\sigma_{x,y}$, $(1,t+1)$ is a transposition in the disjoint cycle decomposition of $\sigma_{y,z}$, and $(2,t+1)$ is a transposition in the disjoint cycle decomposition of $\sigma_{x,z}$. Thus, we may remove each of the transpositions from the disjoint cycle decomposition of the corresponding derangement, and interpret the result as a permutation of $[3,t]$ in each case. Let $\sigma_1$ be the result of removing $(1,2)$ from $\sigma_{x,y}$, $\sigma_2$ be the result of removing $(1,t+1)$ from $\sigma_{yz}$, and $\sigma_3$ be the result of removing $(2,t+1)$ from $\sigma_{x,z}$. Figure~\ref{tetrahedron removal} shows the subgraphs $H_{x,y}, H_{y,z},$ and $H_{x,z}$, from which $\sigma_{x,y}, \sigma_{y,z}, \sigma_{x,z}$ can be read-off, and the constructed permutations $\sigma_1,\sigma_2$ and $\sigma_{3}$ are shown below the corresponding subgraphs.

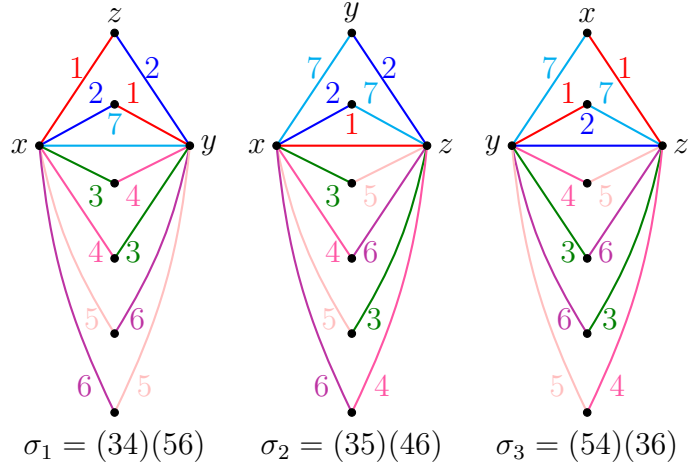
\begin{figure}[h]
 \begin{center}
 \begin{tikzpicture}
 \draw[thick, red] (0,0) -- (1,1.5) node[pos=0.5,above]{$1$};
 \draw[thick, blue] (0,0) -- (1,0.55) node[pos=0.75, above]{$2$};
 \draw[thick, forest] (0,0) -- (1,-0.5) node[pos=0.75,below]{$3$};
 \draw[thick, brilliantrose] (0,0) -- (1,-1.5) node[pos=0.75,below]{$4$};
 \draw[thick, pink] (0,0) to[bend right = 10] (1,-2.5); 
 \draw[pink] (0.7,-2.3) node{$5$};
 \draw[thick, byzantine] (0,0) to[bend right = 10] (1,-3.55);
 \draw[byzantine] (0.6,-3.2) node{$6$};
 
 \draw[thick, blue] (2,0) -- (1,1.5) node[pos=0.5,above]{$2$};
 \draw[thick, red] (2,0) -- (1,0.55) node[pos=0.75,above]{$1$};
 \draw[thick,brilliantrose] (2,0) -- (1,-0.5) node[pos=0.75,below]{$4$};
 \draw[thick, forest] (2,0) -- (1,-1.5) node[pos=0.75,below]{$3$};
 \draw[thick, byzantine] (2,0) to[bend left = 10](1,-2.5); 
 \draw[byzantine] (1.3,-2.3)  node{$6$};
 \draw[thick, pink] (2,0) to[bend left = 10] (1,-3.55);
 \draw[pink] (1.4,-3.2) node{$5$};
 
 \draw[thick, cyan] (0,0) -- (2,0) node[pos=0.5,above]{$7$};
 
 \filldraw (0,0) circle (0.05 cm) node[left]{$x$};
 \filldraw (2,0) circle (0.05 cm) node[right]{$y$};
 \filldraw (1,1.5) circle (0.05 cm) node[above]{$z$};
 \filldraw (1,0.55) circle (0.05 cm);
 
 \filldraw (1,-0.5) circle (0.05 cm);
 \filldraw (1,-1.5) circle (0.05 cm);
 \filldraw (1,-2.5) circle (0.05 cm);
 \filldraw (1,-3.55) circle (0.05 cm);

 \draw (1,-4) node{$\sigma_1 = (3 4) (5 6)$};
 \end{tikzpicture}
 \begin{tikzpicture}
 \draw[thick, cyan] (0,0) -- (1,1.5) node[pos=0.5,above]{$7$};
 \draw[thick, blue] (0,0) -- (1,0.55) node[pos=0.75, above]{$2$};
 \draw[thick, forest] (0,0) -- (1,-0.5) node[pos=0.75,below]{$3$};
 \draw[thick, brilliantrose] (0,0) -- (1,-1.5) node[pos=0.75,below]{$4$};
 \draw[thick, pink] (0,0) to[bend right = 10] (1,-2.5);
 \draw[pink] (0.7, -2.3) node{$5$};
 \draw[thick, byzantine] (0,0) to[bend right = 10]  (1,-3.55);
 \draw[byzantine] (0.6,-3.2) node{$6$};
 
 \draw[thick, blue] (2,0) -- (1,1.5) node[pos=0.5,above]{$2$};
 \draw[thick, cyan] (2,0) -- (1,0.55) node[pos=0.75,above]{$7$};
 \draw[thick,pink] (2,0) -- (1,-0.5) node[pos=0.75,below]{$5$};
 \draw[thick, byzantine] (2,0) -- (1,-1.5) node[pos=0.75,below]{$6$};
 \draw[thick, forest] (2,0) to[bend left = 10]  (1,-2.5);
 \draw[forest] (1.3, -2.3) node{$3$};
 \draw[thick, brilliantrose] (2,0) to[bend left = 10] (1,-3.55);
 \draw[brilliantrose] (1.4,-3.2) node{$4$};
 
 \draw[thick, red] (0,0) -- (2,0) node[pos=0.5,above]{$1$};
 
 \filldraw (0,0) circle (0.05 cm) node[left]{$x$};
 \filldraw (2,0) circle (0.05 cm) node[right]{$z$};
 \filldraw (1,1.5) circle (0.05 cm) node[above]{$y$};
 \filldraw (1,0.55) circle (0.05 cm);
 
 \filldraw (1,-0.5) circle (0.05 cm);
 \filldraw (1,-1.5) circle (0.05 cm);
 \filldraw (1,-2.5) circle (0.05 cm);
 \filldraw (1,-3.55) circle (0.05 cm);
 \draw (1,-4) node{$\sigma_2 = (3 5) (4 6)$};
 \end{tikzpicture}
 \begin{tikzpicture}
  \draw[thick, cyan] (0,0) -- (1,1.5) node[pos=0.5,above]{$7$};
 \draw[thick, red] (0,0) -- (1,0.55) node[pos=0.75, above]{$1$};
 \draw[thick, brilliantrose] (0,0) -- (1,-0.5) node[pos=0.75,below]{$4$};
 \draw[thick, forest] (0,0) -- (1,-1.5) node[pos=0.75,below]{$3$};
 \draw[thick, byzantine] (0,0) to[bend right = 10]  (1,-2.5);
 \draw[byzantine] (0.7,-2.3) node{$6$};
 \draw[thick, pink] (0,0) to[bend right = 10]  (1,-3.55); 
 \draw[pink] (0.6,-3.2) node{$5$};
 
 \draw[thick, red] (2,0) -- (1,1.5) node[pos=0.5,above]{$1$};
 \draw[thick, cyan] (2,0) -- (1,0.55) node[pos=0.75,above]{$7$};
 \draw[thick,pink] (2,0) -- (1,-0.5) node[pos=0.75,below]{$5$};
 \draw[thick, byzantine] (2,0) -- (1,-1.5) node[pos=0.75,below]{$6$};
 \draw[thick, forest] (2,0) to[bend left = 10] (1,-2.5);
 \draw[forest] (1.3,-2.3) node{$3$};
 \draw[thick, brilliantrose] (2,0) to[bend left = 10]  (1,-3.55);
 \draw[brilliantrose] (1.4,-3.2) node{$4$};
 
 \draw[thick, blue] (0,0) -- (2,0) node[pos=0.5,above]{$2$};
 
 \filldraw (0,0) circle (0.05 cm) node[left]{$y$};
 \filldraw (2,0) circle (0.05 cm) node[right]{$z$};
 \filldraw (1,1.5) circle (0.05 cm) node[above]{$x$};
 \filldraw (1,0.55) circle (0.05 cm);
 
 \filldraw (1,-0.5) circle (0.05 cm);
 \filldraw (1,-1.5) circle (0.05 cm);
 \filldraw (1,-2.5) circle (0.05 cm);
 \filldraw (1,-3.55) circle (0.05 cm);

 \draw (1,-4) node{$\sigma_3 = (5 4) (3 6)$};
 \end{tikzpicture}
 \caption{The construction of permutations $\sigma_1,\sigma_2,$ and $\sigma_3$}\label{tetrahedron removal}
 \end{center}
 \end{figure}
 
The resulting permutations $\sigma_1,\sigma_2$ and $\sigma_3$ all remain derangements (now on the set $[3,t]$) and all remain products of disjoint transpositions. Furthermore, since they are on the same set, they are composable. To understand the result of the compositions, we interpret these permutations combinatorially. Letting $U = (N_H(x) \cap N_H(y) \cap N_H(z))\setminus \{v\}$, the definition of $\sigma_1$ implies that it is the restriction of $\sigma_{x,y}$ to the set $c(\{xu : u \in U\})$, having as its image the set $c(\{yu:u\in U\})$, and likewise for $\sigma_2$ and $\sigma_3$. If we consider the value of $\sigma_2 (\sigma_1(a))$ for given color $a \in [3,t]$, we see that there is a vertex $u \in U$ for which $c(xu) = a$ and $c(yu) = \sigma_1(a)$. The color of the edge $zu$ can be described as both $\sigma_2(\sigma_1(a))$ and $\sigma_3(a)$. Hence, the two are equal and we have $\sigma_2 \sigma_1 = \sigma_3$. Each $\sigma_i$ is self-inverse, being a product of disjoint transpositions. Therefore, $\sigma_2 \sigma_1 = \sigma_3$ implies $\sigma_1 \sigma_2 = \sigma_3$ by taking the inverse of both sides. Thus $\sigma_2 \sigma_1 =\sigma_1 \sigma_2$ and so $\sigma_2 \sigma_1 \sigma_2^{-1} = \sigma_1$. 

This equation has an important implication for the disjoint cycle decomposition of $\sigma_1$. Let $(a_1, b_1),\cdots ,(a_k, b_k)$ be the transpositions in the disjoint cycle decomposition of $\sigma_1$, so that $\sigma_1 = (a_1, b_1) \cdots (a_k, b_k)$. Note that we have $k = \frac{t-2}{2}$, since every element of $[3,t]$ appears in exactly one of the transpositions.  Then $\sigma_2 \sigma_1 \sigma_2^{-1} = \sigma_1$ implies 
$$
(\sigma_2(a_1), \sigma_2(b_1))\cdots (\sigma_2(a_k), \sigma_2(b_k)) = (a_1, b_1) \cdots (a_k, b_k).
$$
Letting $T = \{ (a_1, b_1), \cdots , (a_k, b_k) \}$, the above equation implies that $(a_1, b_1) \mapsto (\sigma_2(a_1), \sigma_2(b_1))$ defines a permutation of $T$, say $\mu$. This permutation $\mu$ is also a derangement of $T$, and is also a product of disjoint transpositions. To show the first claim, suppose $\mu(a_i, b_i) = (a_i, b_i)$. Then we have, by definition of $\mu$, either $\sigma_2(a_i) = a_i$ or $\sigma_2(a_i) = b_i$. In the first case, $\sigma_2$ would not be a derangement. In the second case $\sigma_1 (\sigma_2(a_i)) = \sigma_1(b_i) = a_i$, as $(a_i, b_i)$ is one of the cycles in the disjoint cycle decomposition of $\sigma_1$. But then $\sigma_3(a_i) = a_i$, as $\sigma_1 \sigma_2 = \sigma_3$. This contradicts the fact that $\sigma_3$ is a derangement. Thus, neither case is possible and $\mu$ is a derangement of $T$. To show that $\mu$ is a product of disjoint transpositions, it suffices to show that $\mu^2$ is the identity. This follows from the observation that $\mu^2 (a_i, b_i) = \mu( \sigma_2(a_i), \sigma_2(b_i)) = (\sigma_2^2(a_i), \sigma_2^2(b_i)) = (a_i, b_i)$. Hence $\mu$ is a derangement of $T$ and is a product of disjoint transpositions. Therefore, the size of $T$ is even. That is, $\frac{t-2}{2}$ is even, which implies $t \equiv 2 \pmod{4}$.
\end{proof}
 
We will now prove, in stages, that if $G$ has either three vertices of degree $t+1$ or two vertices of degree $t+1$ and a vertex of degree $t$, then $t \equiv 2 \pmod{4}$. This implies Lemma~\ref{structure}. To do this, we will produce a graph $H$ with a good coloring from $G$ and $c$, and which satisfies the hypotheses of Lemma~\ref{high degree in good}. In the case of three vertices of degree $t+1$, the graph $H$ will be a subgraph of $G$. In the other case, $H$ will not be a subgraph of $G$ but will instead be constructed from one. Showing that this construction satisfies the necessary conditions will be the bulk of the proof.

In light of this, we may assume that $G$ has at least two vertices of degree $t+1$. We fix two of these and label them $x$ and $y$. From the remaining vertices in $G$, choose $z$ to have maximum degree. The case of three vertices of degree $t+1$ is then precisely the case where $z$ has degree $t+1$, and the case of two vertices of degree $t+1$ and a vertex of degree $t$ is precisely the case where $z$ has degree $t$.

The first case is easy:

\begin{lemma}\label{z deg t + 1}
Suppose $z$ has degree $t+1$ in $G$. Then $t \equiv 2 \pmod{4}$. 
\end{lemma}
\begin{proof}
Let $H = (V(G), E(G_{x,y}) \cup E(G_{y,z}) \cup E(G_{x,z}))$, noting that each of $G_{x,y}, G_{y,z},$ and $G_{x,z}$ is a subgraph of $G$ by definition. Any $C_4$ in $H$ contains at least two of the vertices $x,y$ and $z$, as the remaining vertices form an independent set in $H$. Thus, any $C_4$ is contained in one of $G_{x,y}$, $G_{y,z}, G_{x,z}$. Since each of $G_{x,y}$, $G_{y,z}$, and $G_{x,z}$ is good under $c$ (by Lemma~\ref{good restriction}), the graph $H$ contains no trichromatic $C_4$'s under $c$, and the colors seen by the edges incident to $x$ are the same as the colors seen by the edges incident to $y$ and $z$ (by the argument that $C_u=C_v$ in Lemma~\ref{derangements}). Hence, $c$ restricts to a good coloring of $H$ and $H$ has three vertices of degree $t+1$. By Lemma~\ref{high degree in good}, $t \equiv 2 \pmod{4}.$ 
\end{proof}
The case where $z$ has degree $t$ is more complicated. We start with the graph $G' =  (V(G), E(G_{x,y}) \cup E(G_{y,z}) \cup E(G_{x,z}))$, as before. We cannot pass this to Lemma~\ref{high degree in good} yet, as $z$ is not of degree $t+1$. To fix this, we note that there is a unique vertex $u$ which is not adjacent to $z$. We would like to take $H = G' + zu$, but we need to make sure we can extend the coloring $c$ to the edge $zu$. We will do this in two cases. 

Note that, as in the proof of Lemma~\ref{z deg t + 1}, $c$ is a good coloring of $G'$, or else one can find a rainbow $B_{t,3}$-copy in $G$ under $c$. Since $c$ is a good coloring of $G'$, we may assume $c$ uses the colors $[t+1]$, that $c(xz) = 1, c(yz) = 2$ and that $c(xy) = t+1$. We define the vertex $v$ to be the unique vertex satisfying $c(xv) = 2$. We first consider the case where $u = v$, that is, $z$ is not adjacent to the unique vertex $v$ satisfying $c(xv) = 2$. 

\begin{lemma}\label{u=v case}
Suppose that $zv$ is not an edge of $G$ and that $z$ has degree $t$. Then $t \equiv 2 \pmod{4}$.
\end{lemma}
\begin{proof}
We want to set $H = G' + zv$ and extend $c$ to the edge $zv$, but we first need to know what color $c(zv)$ could potentially be. The edges incident to $z$ use $t$ colors, so we want to use the remaining color in $[t+1]$ for $c(zv)$. We claim that this remaining color is $t+1$. Indeed, suppose that $z$ is incident to an edge of color $t+1$, say $zw$. Then we have the cycle $xyzwx$ in $G'$. The first edge in this cycle has color $t+1$, the second has color $2$, and the third has color $t+1$. By the goodness of $c$ on $H$, this cycle must be bichromatic and thus we must have $c(xw) = 2$. This is the definition of $v$, and so $w = v$ and $zv$ would be an edge of $G$. Hence, since $zv$ is not an edge, the color in $[t+1]$ that $z$ is not incident to must be $t+1$.

Noting that $v$ is not incident to an edge of color $t+1$, we can now extend $c$ to $G'+ zv$ by setting $c(zv) = t+1$ and obtain a proper edge-coloring of $G'+zv$. We need to show that this extension is a good coloring of $G' + zv$.

\begin{figure}[h!]
\begin{center}
\begin{tikzpicture}
\draw[thick,red] (0,0) -- (1,1) node[pos=0.5, left]{\scriptsize$1$};
\draw[thick,blue] (1,1) -- (2,0) node[pos=0.5, right]{\scriptsize$t+1$};
\draw[thick,red] (2,0) -- (1,-1) node[pos=0.5, right]{\scriptsize$1$};
\draw[thick,blue] (1,-1) -- (0,0) node[pos=0.5, left]{\scriptsize$t+1$};

\filldraw (0,0) circle (0.05 cm) node[left] {$x$};
\filldraw (1,1) circle (0.05 cm) node[above] {$z$};
\filldraw (2,0) circle (0.05 cm) node[right] {$v$};
\filldraw (1,-1) circle (0.05 cm) node[below]{$y$};
\end{tikzpicture}
\begin{tikzpicture}
\draw[thick,forest] (0,0) -- (1,1) node[pos=0.5, left]{\scriptsize$2$};
\draw[thick,blue] (1,1) -- (2,0) node[pos=0.5, right]{\scriptsize$t+1$};
\draw[thick,forest] (2,0) -- (1,-1) node[pos=0.5, right]{\scriptsize$2$};
\draw[thick,blue] (1,-1) -- (0,0) node[pos=0.5, left]{\scriptsize$t+1$};

\filldraw (0,0) circle (0.05 cm) node[left] {$y$};
\filldraw (1,1) circle (0.05 cm) node[above] {$z$};
\filldraw (2,0) circle (0.05 cm) node[right] {$v$};
\filldraw (1,-1) circle (0.05 cm) node[below]{$x$};
\end{tikzpicture}
\begin{tikzpicture}[scale=1.2]
\draw[thick, brilliantrose] (0,0) -- (1,1) node[pos=0.5, left]{\scriptsize$a$};
\draw[thick, blue] (1,1) -- (2,0) node[pos=0.5, right]{\scriptsize$t+1$};
\draw[thick, forest] (2,0) -- (1,-1) node[pos=0.5, right]{\scriptsize$2$};
\draw[thick, pink] (1,-1) -- (0,0) node[pos=0.5, left]{\scriptsize$b$};
\draw[thick, forest] (1,1) -- (1,0) node[pos=0.75, right]{\scriptsize$2$};
\draw[thick, blue] (1,-1) -- (1,0) node[pos=0.75, right]{\scriptsize $t+1$};

\filldraw (0,0) circle (0.05 cm);
\filldraw (1,1) circle (0.05 cm) node[above]{$z$};
\filldraw (2,0) circle (0.05 cm) node[right]{$v$};
\filldraw (1,-1) circle (0.05 cm) node[below]{$x$};

\filldraw (1,0) circle (0.05 cm) node[above left]{$y$};
\end{tikzpicture}
\begin{tikzpicture}[scale=1.2]
\draw[thick, brilliantrose] (0,0) -- (1,1) node[pos=0.5, left]{\scriptsize$a$};
\draw[thick, blue] (1,1) -- (2,0) node[pos=0.5, right]{\scriptsize$t+1$};
\draw[thick, red] (2,0) -- (1,-1) node[pos=0.5, right]{\scriptsize$1$};
\draw[thick, pink] (1,-1) -- (0,0) node[pos=0.5, left]{\scriptsize$b$};
\draw[thick, red] (1,1) -- (1,0) node[pos=0.75, right]{\scriptsize$1$};
\draw[thick, blue] (1,-1) -- (1,0) node[pos=0.75, right]{\scriptsize $t+1$};

\filldraw (1,0) circle (0.05 cm) node[above left]{$x$};
\filldraw (0,0) circle (0.05 cm);
\filldraw (1,1) circle (0.05 cm) node[above]{$z$};
\filldraw (2,0) circle (0.05 cm) node[right]{$v$};
\filldraw (1,-1) circle (0.05 cm) node[below]{$y$};
\end{tikzpicture}

\caption{$C_4$-copies in $G' + zv$ containing $zv$}\label{z C4 copies}
\end{center}
\end{figure}
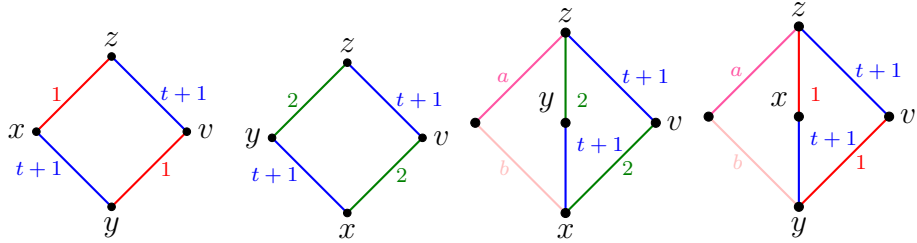

To do this, we need to show that $G' + zv$ contains no trichromatic four-cycle. Four-cycles in $G' + zv$  that don't contain $zv$ are also four-cycles in $G'$, and hence are not trichromatic. So we consider a four-cycle containing $zv$. The possible forms of such a cycle are depicted in Figure~\ref{z C4 copies}. 

In the first two cases, the colors are determined and the cycles are bichromatic. In the last two cases, the cycles are rainbow by the properness of $c$ on $G' + zv$. So $c$ is good on $G' +zv$. With $H = G' + zv$, we now satisfy the hypotheses of Lemma~\ref{high degree in good}, and thus $t \equiv 2 \pmod{4}$.
\end{proof}

Now we prove the full thing.
\begin{lemma}
If $z$ has degree $t$, then $t \equiv 2 \pmod{4}$.
\end{lemma}
\begin{proof}
By Lemma~\ref{u=v case}, we are done if $zv$ is not an edge in $G$, so we will assume that it is. As in the proof of Lemma~\ref{high degree in good}, we know $c(zv) = t+1$. Again, we let $u$ be the unique vertex in $G'$ not incident to  $z$. 

Our goal is to extend the coloring $c$ to edge $zu$ so that we may apply Lemma~\ref{high degree in good} to $G' + zu$. We will not add this edge yet, and instead work with within $G'$ for the first part of the proof. 

Note that $t$ of the colors in $[t+1]$ appear on an edge incident to $z$. We want to extend $c$ by setting $c(zu)$ to be the unique color in $[t+1]$ which is not seen by $z$ in $G'$.     By permuting $\{3, \dots, t\}$ if necessary, we may assume $c(xu) = 3$ and $c(yu) = 4$. In order for the extension to be proper, we will need to know that the missing color at $z$ is neither $3$ nor $4$.

We will now show that $z$ is incident to edges of colors $3$ and $4$, which will ensure that we get a proper coloring when we extend $c$.
Suppose that $z$ is not incident to an edge of the color 3. Since this is the only color not seen by $z$, there is some vertex $w$ such that $c(wz) = 4$. Let $c_1 = c(yw)$. 

If $c_1 \neq 3$, then $c_1 = c(zs)$ for some vertex $s$. Applying the goodness of $c$ on $G'$ to the cycle $yszwy,$ we have $c(ys) = 4$. But $s$ is not $u$, as $z$ is adjacent to one and not the other. Hence, $ys$ and $yu$ are distinct edges of the same color, contradicting the properness of $c$.
\begin{figure}[h]
\begin{center}
 \begin{tikzpicture}
 \draw[thick, red] (0,0) -- (1,1) node[pos=0.5, above]{$c_1$}; 
 \draw[thick] (0,0) -- (1,-1); 
 \draw[thick, blue] (2,0) -- (1,1) node[pos=0.5, above]{$4$}; 
 \draw[thick, red] (2,0) -- (1,-1) node[pos=0.5, below]{$c_1$};
 \draw[thick, blue] (0,0) -- (1,0) node[pos=0.75, above]{$4$};
 
 \filldraw (0,0) circle (0.05 cm) node[left]{$y$};
 \filldraw (2,0) circle (0.05 cm) node[right]{$z$};
 \filldraw (1,1) circle (0.05 cm) node[above]{$w$};
 \filldraw (1,-1) circle (0.05 cm) node[below]{$s$};
 \filldraw (1,0) circle (0.05 cm) node[right]{$u$};
 \end{tikzpicture}
\end{center}
\caption{Contradiction proving $c_1 = 3$}
\end{figure}
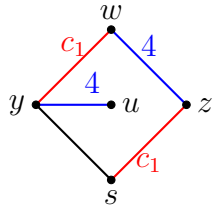

Hence, we must have $c(yw) = c_1 = 3$. Applying the goodness of $c$ on $G'$ to the cycle $xuywx$ yields $c(xw) = 4$. But then $c(xw) = c(zw) = 4$, contradicting properness. Hence, $z$ must be incident to an edge of color $3$. By a symmetric argument with $y$, $z$ must be incident to an edge of color $4$. At this point, we know that the desired extension of $c$ will be proper.
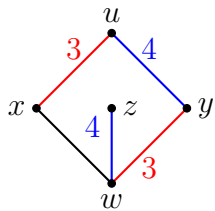
\begin{figure}[h]
\begin{center}
 \begin{tikzpicture}
 \draw[thick, red] (0,0) -- (1,1) node[pos=0.5, above]{$3$}; 
 \draw[thick] (0,0) -- (1,-1); 
 \draw[thick, blue] (2,0) -- (1,1) node[pos=0.5, above]{$4$}; 
 \draw[thick, red] (2,0) -- (1,-1) node[pos=0.5, below]{$3$};
 \draw[thick, blue] (1,-1) -- (1,0) node[pos=0.75, left]{$4$};
 
 \filldraw (0,0) circle (0.05 cm) node[left]{$x$};
 \filldraw (2,0) circle (0.05 cm) node[right]{$y$};
 \filldraw (1,1) circle (0.05 cm) node[above]{$u$};
 \filldraw (1,-1) circle (0.05 cm) node[below]{$w$};
 \filldraw (1,0) circle (0.05 cm) node[right]{$z$};
 \end{tikzpicture}
\end{center}
\caption{Contradiction resulting from $c_1 = 3$}
\end{figure}

We can now describe the missing color.
Let $r$ be the vertex for which $c(zr) = 4$ and let $q$ be the vertex for which $c(zq) = 3$. Let $c_2 = c(xq)$ and $c_3 = c(yr)$. Note that $c_2 \neq 3$ and $c_3 \neq 4$ by the properness of $c$.

If $c(zp) = c_2$ for some vertex $p$, then applying the goodness of $c$ on $G'$ to the cycle $xqzpx$ yields $c(xp) = 3$, since $c(xq) = c_2, c(qz) = 3,$ and $c(zp) = c_2$. 

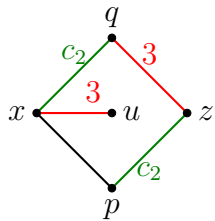
\begin{figure}[h]
\begin{center}
 \begin{tikzpicture}
 \draw[thick, forest] (0,0) -- (1,1) node[pos=0.5, above]{$c_2$}; 
 \draw[thick] (0,0) -- (1,-1); 
 \draw[thick, red] (2,0) -- (1,1) node[pos=0.5, above]{$3$}; 
 \draw[thick, forest] (2,0) -- (1,-1) node[pos=0.5, below]{$c_2$};
 \draw[thick, red] (0,0) -- (1,0) node[pos=0.75, above]{$3$};
 
 \filldraw (0,0) circle (0.05 cm) node[left]{$x$};
 \filldraw (2,0) circle (0.05 cm) node[right]{$z$};
 \filldraw (1,1) circle (0.05 cm) node[above]{$q$};
 \filldraw (1,-1) circle (0.05 cm) node[below]{$p$};
 \filldraw (1,0) circle (0.05 cm) node[right]{$u$};
 \end{tikzpicture}
\end{center}
\caption{Contradiction proving the missing color is $c_2$}
\end{figure}

But then $p = u$ by properness, contradicting the fact that $z$ and $u$ are not adjacent. So $z$ does not see an edge of color $c_2$. A similar argument, using $y$ and $r$ instead of $x$ and $q$, shows that $z$ does not see an edge of color $c_3$. Thus $c_2 = c_3$, and this is the missing color. 

Let $H = G' + zu$ and extend $c$ to $H$ by defining $c(zu) = c_2$. We have already noted that this is a proper edge-coloring, and it still uses at most $t+1$ colors. To show $c$ is good, we need to show that there are no trichromatic four-cycles in $H$ under $c$.

\begin{figure}[h]
\begin{center}
\begin{tikzpicture}
\draw[thick,byzantine] (0,0) -- (1,1) node[pos=0.5, left]{$t+1$};
\draw[thick,cyan] (1,1) -- (2,0) node[pos=0.5, right]{$2$};
\draw[thick,forest] (2,0) -- (1,-1) node[pos=0.5, right]{$c_2$};
\draw[thick,red] (1,-1) -- (0,0) node[pos=0.5, left]{$3$};

\filldraw (0,0) circle (0.05 cm) node[left] {$x$};
\filldraw (1,1) circle (0.05 cm) node[above] {$y$};
\filldraw (2,0) circle (0.05 cm) node[right] {$z$};
\filldraw (1,-1) circle (0.05 cm) node[below]{$u$};
\end{tikzpicture}
\begin{tikzpicture}
\draw[thick,byzantine] (0,0) -- (1,1) node[pos=0.5, left]{$t+1$};
\draw[thick,pink] (1,1) -- (2,0) node[pos=0.5, right]{$1$};
\draw[thick,forest] (2,0) -- (1,-1) node[pos=0.5, right]{$c_2$};
\draw[thick,blue] (1,-1) -- (0,0) node[pos=0.5, left]{$4$};

\filldraw (0,0) circle (0.05 cm) node[left] {$y$};
\filldraw (1,1) circle (0.05 cm) node[above] {$x$};
\filldraw (2,0) circle (0.05 cm) node[right] {$z$};
\filldraw (1,-1) circle (0.05 cm) node[below]{$u$};
\end{tikzpicture}
\begin{tikzpicture}
\draw[thick] (0,0) -- (1,1);
\draw[thick] (1,1) -- (2,0);
\draw[thick, forest] (2,0) -- (1,-1) node[pos=0.5, right]{$c_2$};
\draw[thick, red] (1,-1) -- (0,0) node[pos=0.5, left]{$3$};
\filldraw (0,0) circle (0.05 cm) node[left]{$x$};
\filldraw (1,1) circle (0.05 cm) node[above]{$s$};
\filldraw (2,0) circle (0.05 cm) node[right]{$z$};
\filldraw (1,-1) circle (0.05 cm) node[below]{$u$};
\end{tikzpicture}
\begin{tikzpicture}
\draw[thick] (0,0) -- (1,1);
\draw[thick] (1,1) -- (2,0);
\draw[thick, forest] (2,0) -- (1,-1) node[pos=0.5, right]{$c_2$};
\draw[thick, blue] (1,-1) -- (0,0) node[pos=0.5, left]{$4$};
\filldraw (0,0) circle (0.05 cm) node[left]{$y$};
\filldraw (1,1) circle (0.05 cm) node[above]{$s$};
\filldraw (2,0) circle (0.05 cm) node[right]{$z$};
\filldraw (1,-1) circle (0.05 cm) node[below]{$u$};
\end{tikzpicture}
\caption{The four kinds of $C_4$'s containing $zv$, where $s \not\in \{x,y\}$}
\end{center}
\end{figure}
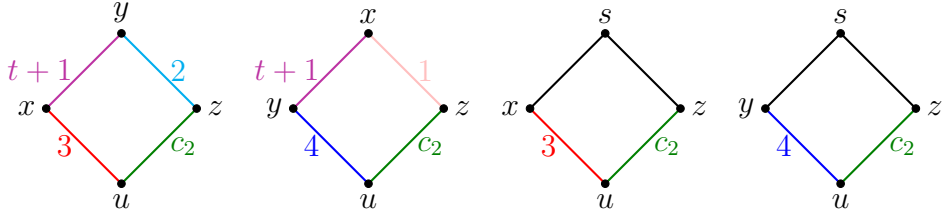

Any four-cycle not containing the edge $zu$ is also a four-cycle in $G$. By the goodness of $c$ on $H$, none of these are trichromatic. Any four-cycle containing $zu$ must also contain one of the vertices $x$ and $y$, as the only neighbors of $u$ in $H$ are $x,y$ and $z$. These four-cycles are $xyzux, xyuzx$, and the cycles of the forms $xszux$ and $yszuy$ where $s$ is a vertex other than $x$ and $y$.
 
Consider the four-cycles $xyzux$ and $xyuzx$. Since $u$ is not $v$ and $c(zv) = t+1$, we know $c(zu) = c_2 \neq t+1$. Hence, these two cycles are rainbow. 

Now consider the cycles of the form $xszux$ where $s$ is not $y$. As $c(xu) = 3$ and $c(zu) = c_2$, such a cycle is trichromatic only if $c(zs) = 3$ or $c(xs) = c_2$. But recall that the vertex $q$ was defined by $c(zq) = 3$, and $c_2$ was defined by $c_2 = c(xq)$. So in either case, $s = q$. The cycle $xqzux$ is bichromatic, and so no cycle of the form $xszux$ is trichromatic. Similarly, if a cycle of the form $yszuy$ is trichromatic, it must be $yrzuy$ in particular. This cycle is bichromatic, and so $c$ is a good coloring of $H$.

The graph $H$ with good coloring $c$ now satisfies the hypotheses of Lemma~\ref{high degree in good}, and so $t \equiv 2 \pmod{4}$
\end{proof}
Lemma~\ref{structure} now follows. If $G$ has more than two vertices of degree $t+1$, or two vertices of degree $t+1$ and a vertex of degree $t$, then $t \equiv 2 \pmod{4}.$ Hence if $t \equiv 0 \pmod{4}$, the conclusion of Lemma~\ref{structure} follows.

\subsection{Lower bound when $t = 2^s$}\label{power of 2 construction}

\begin{lemma}
    Let $s \geq 2$ and suppose $t = 2^s$. Then, $\mathrm{ex^*}(n,B_{t,3}) \geq \frac{t}{2}n + O(1)$.
\end{lemma}
\begin{proof}
    We will show that there exists a proper edge-coloring of $G:= K_{t,t}$ containing no rainbow $B_{t,3}$. Let $X = \{x_1, \dots, x_t\}, Y = \{y_1, \dots, y_t\}$ be the bipartition of $G$ and identify each of $X$ and $Y$ with a copy of $\mathbb F_2^s$. Let us define an edge coloring $c:E(G) \to \mathbb{F}_2^s$ by $c(x_iy_j) = x_i - y_j$. Clearly, this is a proper edge coloring. Suppose there exists some rainbow $B_{t,3}$, say $B$ in $G$. Without loss of generality, we may assume the beginning of the path of $B$ starts in $X$. After relabeling, we may assume the path $x_1,y_1,x_2,y_2$ is the handle of $B$ and the vertices $x_3, \dots, x_{t-1}$ are incident to the bristles of $B$. Observe that $\{c(y_2x_1), c(y_2x_t)\} = \{x_1 - y_1, y_1 - x_2\}$ because $B$ is rainbow and $y_2$ is incident to an edge of each color. Because the coloring is proper, this implies that $c(y_2x_1) = y_1 - x_2$. But then, this means that $y_2 + x_2 = x_1 + y_1$. Since we are in a field of characteristic $2$, we have $y_2 - x_2 = x_1 - y_1$. That is, $c(y_2x_2) = c(x_1y_1)$, a contradiction to $B$ being rainbow. Therefore, $G$ contains no rainbow $B_{t,3}$. \\
    Now for $n \in \mathbb{N}$, taking $\left \lfloor \frac{n}{2t} \right \rfloor$ disjoint copies of $G$ along with some isolated vertices and coloring $G$ as described above, we have
    $$\mathrm{ex^*}(n,B_{t,3}) \geq \left \lfloor \frac{n}{2t} \right \rfloor t^2 = \frac{t}{2}n + O(1), $$
as desired.
\end{proof}

\subsection{Lower bound when $t = 3^s - 1$}\label{power of 3 construction}

In this section we prove Theorem~\ref{power of 3 lower bound}, which we restate here for convenience.
\powerofthreelowerbound*
\begin{proof}
    We will show that there exists a proper edge-coloring of $G:= K_{t + 1}$ containing no rainbow $B_{t,3}$. Set $V(G) = \mathbb{F}_3^s$ and define the edge-coloring $c: E(G) \to \mathbb{F}_3^s$ by $c(uv) = u + v$. Clearly, this is a proper edge coloring. Suppose there exists some rainbow $B_{t,3}$, say $B$ in $G$. Let $zyxv$ be the handle of $B$. Without loss of generality, we may assume that $v = 0$. Indeed, if $v \neq 0$, then we may find another rainbow broom in $G$ with base at $0$ by subtracting $v$ from every vertex in $B$. Now, by our choice of coloring, we have $c(vx) = x$, $c(xy) = x + y$ and $c(yz) = y + z$. Furthermore, since $v$ is incident to every color except $0$,
    $$\{ c(vy), c(vz) \} \cap \{c(xy), c(yz)\} = \{y, z\} \cap \{x + y, y + z\} \neq \emptyset.$$
    Now $y\ne x+y,y\ne y+z$, and $z\ne y+z$, so this implies $z=x+y$. Then $y+z$ is the color not incident to $v$, so $y+z=0$. Hence $0=y+z=x+2y=x-y$ which contradicts $x\ne y$.
    Now for $n \in \mathbb{N}$, taking $\left \lfloor \frac{n}{t + 1} \right \rfloor$ disjoint copies of $G$ along with some isolated vertices and coloring $G$ as described above, we have
    $$\mathrm{ex^*}(n,B_{t,3}) \geq \left \lfloor \frac{n}{t + 1} \right \rfloor \binom{t + 1}{2} = \frac{t}{2}n + O(1).$$
\end{proof}

\subsection{Rainbow $B_{10,3}$-free cliques}

So far, we have established that often, $\mathrm{ex}^*(n,B_{t,3}) = \frac{t}{2}n + O(1)$, and almost always, we have found an extremal construction which looks like disjoint copies of $K_{t+1}$ equipped with a ``good'' edge-coloring. While the construction we give to show that $\mathrm{ex}^*(n,B_{t,3}) = \frac{t}{2}n + O(1)$ when $t = 2^s$ is not of this form, it is not obvious that copies of $K_{t+1}$ \textit{cannot} be edge-colored to avoid a rainbow $B_{t,3}$ when $t = 2^s$. (Indeed, when $s = 3$ we have $2^s = 8$, and Theorem~\ref{power of 3 lower bound} shows that $K_9$ \textit{can} be edge-colored without a rainbow $B_{8,3}$.)   

Thus, it might be reasonable to conjecture that $K_{t+1}$ can always be properly edge-colored while avoiding a rainbow $B_{t,3}$-copy. If this were true, then several strong statements would follow: in light of Theorem~\ref{multiple of 4 upper bound}, we would have $\mathrm{ex}^*(n,B_{t,3}) = \frac{t}{2}n + O(1)$ whenever $t \equiv 0 \pmod 4$, and in light of Lemma~\ref{extremal space restriction}, the problem of determining $\mathrm{ex}^*(n,B_{t,3})$ in general would be reduced to determining the densest subgraph of $K_{t+2}$ which can be properly edge-colored without a rainbow $B_{t,3}$-copy. The purpose of this subsection is to demonstrate that the problem cannot be so reduced. Indeed, every proper edge-coloring of $K_{11}$ contains a rainbow $B_{10,3}$-copy.

\begin{prop}\label{B10structure}
No proper edge-coloring of $K_{11}$ is rainbow $B_{10,3}$-free.
\end{prop}
\begin{proof}
Let $c$ be a proper edge-coloring of $K_{11}$. Suppose that there is no rainbow $B_{10,3}$-copy under $c$. Consider a vertex $x$ in $K_{11}$. We may assume that the colors incident to $x$ are $[10]$. Suppose that, for some $wuv$ in $V(K_{11})-x$, $c(wv)$ and $c(uv)$ are not in $[10]$. Then $wuvx$ forms a handle of a rainbow broom whose bristles are the remaining neighbors of $x$. Hence, the edges with colors not in $[10]$ must form a matching, as $x$ is a universal vertex. We may recolor these edges without violating properness so that they all have the same color. If there is a rainbow $B_{10,3}$ after this re-coloring, then there was one before, as the edges which are distinct colors under the re-coloring must be distinct in the original. Hence, we may assume that $c$ uses only the colors from $[11]$, without loss of generality.

Let $M_\ell$ be the matching in color $\ell$, for each $\ell \in [11]$. Note that each matching must have $5$ edges in order for every edge to be colored. Hence, we may uniquely label every vertex by the unique color in $[11]$ it is not incident to. Suppose that we have vertices labeled $A,B,C$ and $D$ and that the edge $BC$ is colored $A$. Then the cycle $ABCD$ cannot be trichromatic. To see this, suppose that it is. Since the vertex labeled $A$ is not incident to color $A$ and since $c$ is proper, the repeated colors must be $c(AB)$ and $c(CD)$. Call their common color $x$ and let $y = c(AD)$. This is shown in Figure~\ref{c4 in k11}. The path $BCDA$ is rainbow since $y \neq A$. Furthermore, of the seven edges from $A$ to a vertex not in $\{B,C,D\}$, none of them are colored $A$, by the definition of the vertex labeling, and none of them are colored $x$ or $y$, as edges of those colors are incident to $A$ in the cycle. Hence, we can form a rainbow $B_{10,3}$-copy. 

\begin{figure}[h]
\begin{center}
\begin{tikzpicture}[
labelstyle/.style={
       circle, draw=black,
       thick, inner sep=2pt, minimum size=1 mm,
       outer sep=0pt
        }, scale=2]
        \node (A) at  (0,0) [labelstyle]{$A$};
        \node (B) at (0,1) [labelstyle]{$B$};
        \node (D) at (1,0) [labelstyle]{$D$};
        \node (C) at (1,1) [labelstyle]{$C$};
        \draw[thick,blue] (A) -- (B) node[left, pos=0.5]{$x$};
        \draw[thick, red] (B) -- (C) node[above,pos=0.5]{$A$};
        \draw[thick,blue] (C) -- (D) node[right,pos=0.5]{$x$};
        \draw[thick,forest] (D) -- (A) node[below, pos=0.5]{$y$};
\end{tikzpicture}
\end{center}
\caption{Forbidden Configuration}\label{c4 in k11}
\end{figure}
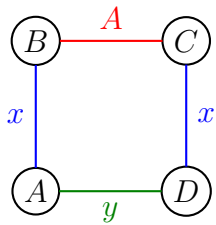

From this observation, we can deduce a useful fact. If $BC$ is an edge of color $A$ and $AB$ is not colored $C$, then we can construct a cycle $ABCD$ where $c(CD) = c(AB)$. The remaining edge, $AD$ cannot be colored $A$ or $c(CD)$, thus giving us a trichromatic cycle. Hence, if $BC$ is an edge of color $A$, then $c(AB) = C$ and $c(AC) = B$. This tells us precisely how the vertex labeled $A$ is incident to the matching $M_A$.

Now we consider the matchings $M_1$ and $M_2$. Up to symmetry, the edges in these two matchings take one of the following forms. $M_1$ is displayed in red and $M_2$ in dashed blue.

\begin{figure}[h]
\begin{center}
\begin{subfigure}{.4\linewidth}
\centering
\begin{tikzpicture}[labelstyle/.style={
       circle, draw=black,
       thick, inner sep=2pt, minimum size=1 mm,
       outer sep=0pt
        }]
       \node (1) at  (0,0) [labelstyle]{$1$};
        \node (2) at (1,0) [labelstyle]{$2$};
        \node (3) at (1,1) [labelstyle]{$3$};
        \node (4) at (2,0) [labelstyle]{$4$};
        \node (9) at (5,0) [labelstyle] {9};
        \node (5) at (2,1) [labelstyle]{5};
        \node (8) at (5,1) [labelstyle] {8};
        \node (6) at (3,2) [labelstyle] {6};
        \node (11) at (3,-1) [labelstyle] {11};
        \node (7) at (4,2) [labelstyle] {7};
        \node (10) at (4,-1) [labelstyle] {10};

        \draw[thick,red] (2) -- (3);
        \draw[thick,blue, dashed] (1) -- (3);

        \draw[thick,red] (4) -- (5);
        \draw[thick,blue, dashed] (5) -- (6);
        \draw[thick,red] (6) -- (7);
        \draw[thick,blue, dashed] (7) -- (8);
        \draw[thick,red] (8) -- (9);
        \draw[thick,blue, dashed] (9) -- (10);
        \draw[thick,red] (10) -- (11);
        \draw[thick,blue, dashed] (11) -- (4);
\end{tikzpicture}
\caption{Case 1}
\end{subfigure}
\begin{subfigure}{.4\linewidth}
\centering
\begin{tikzpicture}[labelstyle/.style={
       circle, draw=black,
       thick, inner sep=2pt, minimum size=1 mm,
       outer sep=0pt
        }]
       \node (1) at  (0,0) [labelstyle]{$1$};
        \node (2) at (1,0) [labelstyle]{$2$};
        \node (3) at (1,1) [labelstyle]{$3$};
        \node (4) at (2,0) [labelstyle]{$4$};
        \node (5) at (2,1) [labelstyle]{5};
        \node (6) at (3,0) [labelstyle]{$6$};
        \node (7) at (3,1) [labelstyle]{$7$};
        \node (8) at (4,0) [labelstyle]{$8$};
        \node (9) at (4,1) [labelstyle]{$9$};
        \node (10) at (5,0) [labelstyle]{$10$};
        \node (11) at (5,1) [labelstyle]{$11$};

        \draw[thick,red] (2) -- (3);
        \draw[thick,red] (4) -- (5);
        \draw[thick,red] (6) -- (7);
        \draw[thick,red] (8) -- (9);
        \draw[thick,red] (10) -- (11);

        \draw[thick,blue, dashed] (1) -- (3);
        \draw[thick,blue, dashed] (4) -- (6);
        \draw[thick,blue, dashed] (5) -- (7);
        \draw[thick,blue, dashed] (8) -- (10);
        \draw[thick,blue, dashed] (9) -- (11);
        \end{tikzpicture}
        \caption{Case 2}
\end{subfigure}
\begin{subfigure}{\linewidth}
\centering
\begin{tikzpicture}[labelstyle/.style={
       circle, draw=black,
       thick, inner sep=2pt, minimum size=1 mm,
       outer sep=0pt
        }]
       \node (1) at  (0,0) [labelstyle]{$1$};
        \node (2) at (1,0) [labelstyle]{$2$};
        \node (3) at (1,1) [labelstyle]{$3$};
        \node (4) at (2,0) [labelstyle]{$4$};
        \node (5) at (2,1) [labelstyle]{5};
        \node (6) at (3,0) [labelstyle]{$6$};
        \node (7) at (3,1) [labelstyle]{$7$};
        \node (8) at (4,0) [labelstyle]{$8$};
        \node (9) at (4,1) [labelstyle]{$9$};
        \node (10) at (5,0) [labelstyle]{$10$};
        \node (11) at (5,1) [labelstyle]{$11$};

        \draw[thick,red] (2) -- (3);
        \draw[thick,red] (4) -- (5);
        \draw[thick,red] (6) -- (7);
        \draw[thick,red] (8) -- (9);
        \draw[thick,red] (10) -- (11);

        \draw[thick,blue, dashed] (1) -- (3);
        \draw[thick,blue, dashed] (5) -- (7);
        \draw[thick,blue, dashed] (6) -- (8);
        \draw[thick,blue, dashed] (9) -- (4);
        \end{tikzpicture}
        \caption{Case 3}
\end{subfigure}
\end{center}
\caption{Possible Relationships between $M_1$ and $M_2$}
\end{figure}
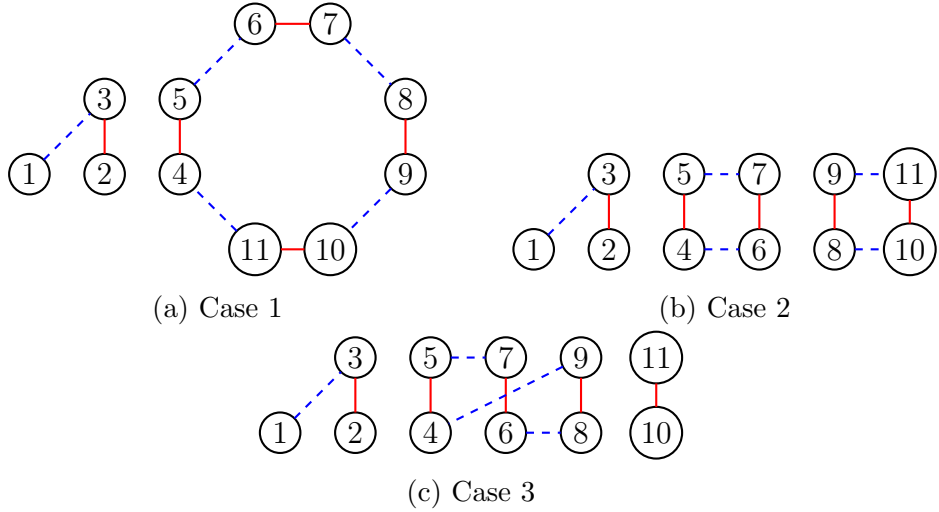

The third case is easily disposed of, as $M_2$ must have five edges in total. The first non-trivial case is Case 1. If we consider the matching $M_3$ in this case, we see that the edge $\{1,2\}$ must have color $3$ by the four-cycle property, and that vertex $3$ is uncovered by $M_3$, by definition. Hence, the remaining edges in $M_3$ must form chords of the octagon in colors $1$ and $2$. These chords must also form a perfect matching of the octagon. There are, up to symmetry, seven ways for a matching of chords to occur. This may be verified by considering the length and relative location of the smallest cycles formed by these chords. The seven types of chordal matchings of an octagon are given in Figure~\ref{perfect chordal matchings}, with $M_3$ in doubled yellow.

\begin{figure}
\begin{subfigure}{0.2\linewidth}
\centering
\begin{tikzpicture}[labelstyle/.style={
       circle, draw=black,
       thick, inner sep=2pt, minimum size=1 mm,
       outer sep=0pt
        }]
        \node (4) at (2,0) [labelstyle]{$4$};
        \node (9) at (5,0) [labelstyle] {9};
        \node (5) at (2,1) [labelstyle]{5};
        \node (8) at (5,1) [labelstyle] {8};
        \node (6) at (3,2) [labelstyle] {6};
        \node (11) at (3,-1) [labelstyle] {11};
        \node (7) at (4,2) [labelstyle] {7};
        \node (10) at (4,-1) [labelstyle] {10};

        \draw[thick,red] (4) -- (5);
        \draw[thick,blue, dashed] (5) -- (6);
        \draw[thick,red] (6) -- (7);
        \draw[thick,blue, dashed] (7) -- (8);
        \draw[thick,red] (8) -- (9);
        \draw[thick,blue, dashed] (9) -- (10);
        \draw[thick,red] (10) -- (11);
        \draw[thick,blue, dashed] (11) -- (4);

        \draw[thick,amber, style=double] (5) -- (7);
        \draw[thick, amber, style=double] (6) -- (8);
        \draw[thick, amber, style=double] (4) -- (10);
        \draw[thick, amber, style=double] (9) -- (11);
\end{tikzpicture}
\caption{Type 1}
\end{subfigure}
\begin{subfigure}{0.2\linewidth}
\centering
\begin{tikzpicture}[labelstyle/.style={
       circle, draw=black,
       thick, inner sep=2pt, minimum size=1 mm,
       outer sep=0pt
        }]
        \node (4) at (2,0) [labelstyle]{$4$};
        \node (9) at (5,0) [labelstyle] {9};
        \node (5) at (2,1) [labelstyle]{5};
        \node (8) at (5,1) [labelstyle] {8};
        \node (6) at (3,2) [labelstyle] {6};
        \node (11) at (3,-1) [labelstyle] {11};
        \node (7) at (4,2) [labelstyle] {7};
        \node (10) at (4,-1) [labelstyle] {10};

        \draw[thick,red] (4) -- (5);
        \draw[thick,blue, dashed] (5) -- (6);
        \draw[thick,red] (6) -- (7);
        \draw[thick,blue, dashed] (7) -- (8);
        \draw[thick,red] (8) -- (9);
        \draw[thick,blue, dashed] (9) -- (10);
        \draw[thick,red] (10) -- (11);
        \draw[thick,blue, dashed] (11) -- (4);

        \draw[thick,amber, style=double] (5) -- (7);
        \draw[thick, amber, style=double] (6) -- (9);
        \draw[thick, amber, style=double] (4) -- (10);
        \draw[thick, amber, style=double] (8) -- (11);
\end{tikzpicture}
\caption{Type 2}
\end{subfigure}
\begin{subfigure}{0.2\linewidth}
\centering
\begin{tikzpicture}[labelstyle/.style={
       circle, draw=black,
       thick, inner sep=2pt, minimum size=1 mm,
       outer sep=0pt
        }]
        \node (4) at (2,0) [labelstyle]{$4$};
        \node (9) at (5,0) [labelstyle] {9};
        \node (5) at (2,1) [labelstyle]{5};
        \node (8) at (5,1) [labelstyle] {8};
        \node (6) at (3,2) [labelstyle] {6};
        \node (11) at (3,-1) [labelstyle] {11};
        \node (7) at (4,2) [labelstyle] {7};
        \node (10) at (4,-1) [labelstyle] {10};

        \draw[thick,red] (4) -- (5);
        \draw[thick,blue, dashed] (5) -- (6);
        \draw[thick,red] (6) -- (7);
        \draw[thick,blue, dashed] (7) -- (8);
        \draw[thick,red] (8) -- (9);
        \draw[thick,blue, dashed] (9) -- (10);
        \draw[thick,red] (10) -- (11);
        \draw[thick,blue, dashed] (11) -- (4);
        
        \draw[thick,amber, style=double] (5) -- (7);
        \draw[thick, amber, style=double] (6) -- (10);
        \draw[thick, amber, style=double] (4) -- (8);
        \draw[thick, amber, style=double] (9) -- (11);
\end{tikzpicture}
\caption{Type 3}
\end{subfigure}
\begin{subfigure}{0.2\linewidth}
\centering
\begin{tikzpicture}[labelstyle/.style={
       circle, draw=black,
       thick, inner sep=2pt, minimum size=1 mm,
       outer sep=0pt
        }]
        \node (4) at (2,0) [labelstyle]{$4$};
        \node (9) at (5,0) [labelstyle] {9};
        \node (5) at (2,1) [labelstyle]{5};
        \node (8) at (5,1) [labelstyle] {8};
        \node (6) at (3,2) [labelstyle] {6};
        \node (11) at (3,-1) [labelstyle] {11};
        \node (7) at (4,2) [labelstyle] {7};
        \node (10) at (4,-1) [labelstyle] {10};

        \draw[thick,red] (4) -- (5);
        \draw[thick,blue, dashed] (5) -- (6);
        \draw[thick,red] (6) -- (7);
        \draw[thick,blue, dashed] (7) -- (8);
        \draw[thick,red] (8) -- (9);
        \draw[thick,blue, dashed] (9) -- (10);
        \draw[thick,red] (10) -- (11);
        \draw[thick,blue, dashed] (11) -- (4);

        \draw[thick,amber, style=double] (5) -- (7);
        \draw[thick, amber, style=double] (6) -- (10);
        \draw[thick, amber, style=double] (4) -- (9);
        \draw[thick, amber, style=double] (8) -- (11);
\end{tikzpicture}
\caption{Type 4}
\end{subfigure}
\begin{subfigure}{0.3\linewidth}
\centering
\begin{tikzpicture}[labelstyle/.style={
       circle, draw=black,
       thick, inner sep=2pt, minimum size=1 mm,
       outer sep=0pt
        }]
        \node (4) at (2,0) [labelstyle]{$4$};
        \node (9) at (5,0) [labelstyle] {9};
        \node (5) at (2,1) [labelstyle]{5};
        \node (8) at (5,1) [labelstyle] {8};
        \node (6) at (3,2) [labelstyle] {6};
        \node (11) at (3,-1) [labelstyle] {11};
        \node (7) at (4,2) [labelstyle] {7};
        \node (10) at (4,-1) [labelstyle] {10};

        \draw[thick,red] (4) -- (5);
        \draw[thick,blue, dashed] (5) -- (6);
        \draw[thick,red] (6) -- (7);
        \draw[thick,blue, dashed] (7) -- (8);
        \draw[thick,red] (8) -- (9);
        \draw[thick,blue, dashed] (9) -- (10);
        \draw[thick,red] (10) -- (11);
        \draw[thick,blue, dashed] (11) -- (4);

        \draw[thick,amber, style=double] (5) -- (8);
        \draw[thick, amber, style=double] (6) -- (11);
        \draw[thick, amber, style=double] (7) -- (10);
        \draw[thick, amber, style=double] (9) -- (4);
\end{tikzpicture}
\caption{Type 5}
\end{subfigure}
\begin{subfigure}{0.3\linewidth}
\centering
\begin{tikzpicture}[labelstyle/.style={
       circle, draw=black,
       thick, inner sep=2pt, minimum size=1 mm,
       outer sep=0pt
        }]
        \node (4) at (2,0) [labelstyle]{$4$};
        \node (9) at (5,0) [labelstyle] {9};
        \node (5) at (2,1) [labelstyle]{5};
        \node (8) at (5,1) [labelstyle] {8};
        \node (6) at (3,2) [labelstyle] {6};
        \node (11) at (3,-1) [labelstyle] {11};
        \node (7) at (4,2) [labelstyle] {7};
        \node (10) at (4,-1) [labelstyle] {10};

        \draw[thick,red] (4) -- (5);
        \draw[thick,blue, dashed] (5) -- (6);
        \draw[thick,red] (6) -- (7);
        \draw[thick,blue, dashed] (7) -- (8);
        \draw[thick,red] (8) -- (9);
        \draw[thick,blue, dashed] (9) -- (10);
        \draw[thick,red] (10) -- (11);
        \draw[thick,blue, dashed] (11) -- (4);

        \draw[thick,amber, style=double] (5) -- (8);
        \draw[thick, amber, style=double] (6) -- (10);
        \draw[thick, amber, style=double] (4) -- (9);
        \draw[thick, amber, style=double] (7) -- (11);
\end{tikzpicture}
\caption{Type 6}
\end{subfigure}
\begin{subfigure}{0.3\linewidth}
\centering
\begin{tikzpicture}[labelstyle/.style={
       circle, draw=black,
       thick, inner sep=2pt, minimum size=1 mm,
       outer sep=0pt
        }]
        \node (4) at (2,0) [labelstyle]{$4$};
        \node (9) at (5,0) [labelstyle] {9};
        \node (5) at (2,1) [labelstyle]{5};
        \node (8) at (5,1) [labelstyle] {8};
        \node (6) at (3,2) [labelstyle] {6};
        \node (11) at (3,-1) [labelstyle] {11};
        \node (7) at (4,2) [labelstyle] {7};
        \node (10) at (4,-1) [labelstyle] {10};

        \draw[thick,red] (4) -- (5);
        \draw[thick,blue, dashed] (5) -- (6);
        \draw[thick,red] (6) -- (7);
        \draw[thick,blue, dashed] (7) -- (8);
        \draw[thick,red] (8) -- (9);
        \draw[thick,blue, dashed] (9) -- (10);
        \draw[thick,red] (10) -- (11);
        \draw[thick,blue, dashed] (11) -- (4);

        \draw[thick,amber, style=double] (5) -- (9);
        \draw[thick, amber, style=double] (6) -- (10);
        \draw[thick, amber, style=double] (7) -- (11);
        \draw[thick, amber, style=double] (4) -- (8);
\end{tikzpicture}
\caption{Type 7}
\end{subfigure}
\caption{Perfect Chordal Matchings of an Octagon, up to symmetry; $M_3$ in doubled yellow} \label{perfect chordal matchings}
\end{figure}
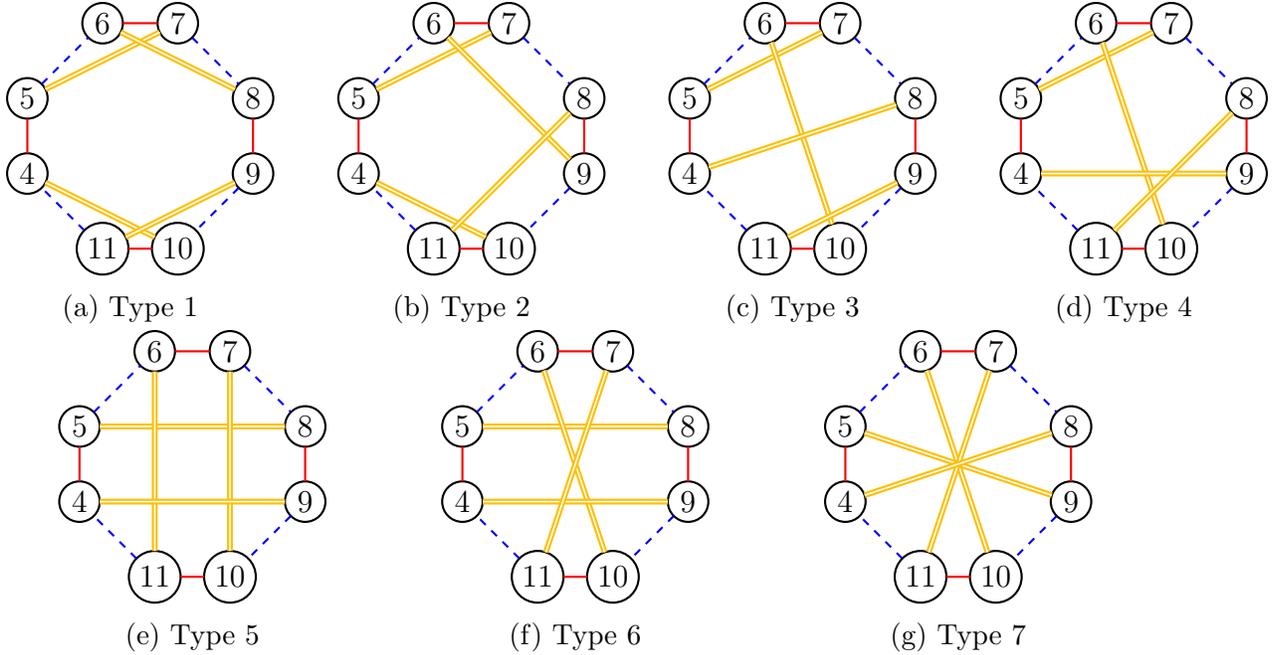

Now we will show that none of these types is possible for $M_3$ in Case 1. 
Suppose that $M_3$ forms a Type 1 perfect chordal matching of the octagon in colors $1$ and $2$. We consider $M_4$. Note that vertex $4$ is not covered by $M_4$, and that the four-cycle fact proved above implies that vertices $1,2,3,5,10$ and $11$ are covered by edges from $M_4$ that are not chords of the octagon. For example, because $c(\{4,5\}) = 1$, we must have $c(\{1,5\}) = 4$. Thus, $M_4$ must contain a perfect matching of the vertices $6,7,8$ and $9$. Observing that $\{6,7\}, \{6,8\}, \{7,8\},$ and $\{8,9\}$ are each contained in one of the matchings $M_1, M_2$ and $M_3$, we see that there are too few remaining edges to form a perfect matching on this set. Hence, $M_4$ cannot form a matching of these vertices, and so $M_3$ cannot form a perfect chordal matching of Type 1. This strategy can be summarized as choosing a color $i$ and then observing that our four-cycle fact shows that $M_i$ must form a perfect matching of exactly four vertices in the octagon. If these vertices cannot admit an additional perfect matching of color $i$, we have a contradiction. Applying this with color $i = 4$ to Types $1$ and $5$, and with color $i = 11$ to Types $2, 3$ and $4$, eliminates these types. We are then left with Types $6$ and $7$. For these types, we repeat the same strategy with color $i = 4$. In these cases, we do not immediately have a contradiction. However, following this strategy shows that there is one possibility for the matching $M_4$ in either of these cases. Once we add this 
matching, we can apply the argument with color $i = 11$ and arrive at a contradiction.

Thus, Case 1 is impossible. 
Note that any two matchings which form a bichromatic $P_4$ must form either Case 1 or Case 3. Both of these lead to contradictions, and thus no bichromatic $P_4$'s may occur under $c$. 

We now turn our attention to Case 2 and consider which edge from $M_3$ covers vertex $7$. Up to symmetry, the two possibilities are $c(\{7,9\}) = 3$ and $c(\{7,4\}) = 3$. Suppose $c(\{7,9\}) = 3$. Vertex $6$ must be covered by an edge of color 3, and the cases up to symmetry are shown in Figure~\ref{73_cases}. Observe that the first three cases contain a bichromatic $P_4$, namely $(6,5,7,9,11)$ in the first, $(5,7,9,11,6)$ in the second, and $(9,7,6,10,11)$ in the third. As such, these cases lead to a contradiction and the only possibility is the fourth case in Figure~\ref{73_cases}. We will return to this case after we have addressed the possibilities when $c(\{7,4\})=3$.

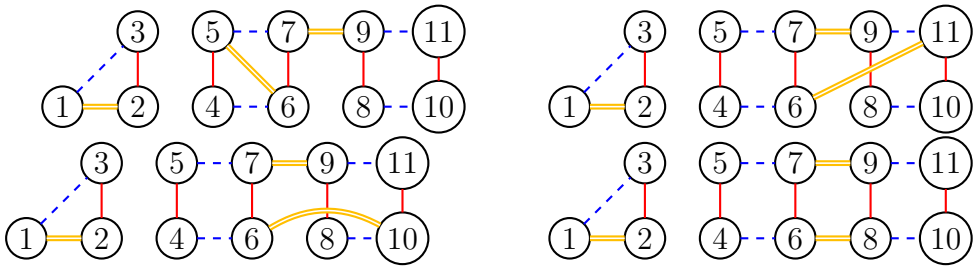
\begin{figure}[h]
\begin{center}
\begin{subfigure}{.4\linewidth}
\centering
\begin{tikzpicture}[labelstyle/.style={
       circle, draw=black,
       thick, inner sep=2pt, minimum size=1 mm,
       outer sep=0pt
        }]
       \node (1) at  (0,0) [labelstyle]{$1$};
        \node (2) at (1,0) [labelstyle]{$2$};
        \node (3) at (1,1) [labelstyle]{$3$};
        \node (4) at (2,0) [labelstyle]{$4$};
        \node (5) at (2,1) [labelstyle]{5};
        \node (6) at (3,0) [labelstyle]{$6$};
        \node (7) at (3,1) [labelstyle]{$7$};
        \node (8) at (4,0) [labelstyle]{$8$};
        \node (9) at (4,1) [labelstyle]{$9$};
        \node (10) at (5,0) [labelstyle]{$10$};
        \node (11) at (5,1) [labelstyle]{$11$};

        \draw[thick,red] (2) -- (3);
        \draw[thick,red] (4) -- (5);
        \draw[thick,red] (6) -- (7);
        \draw[thick,red] (8) -- (9);
        \draw[thick,red] (10) -- (11);

        \draw[thick,blue, dashed] (1) -- (3);
        \draw[thick,blue, dashed] (4) -- (6);
        \draw[thick,blue, dashed] (5) -- (7);
        \draw[thick,blue, dashed] (8) -- (10);
        \draw[thick,blue, dashed] (9) -- (11);

        \draw[thick, amber, style=double] (7) -- (9);
        \draw[thick,amber, style=double] (6) -- (5);
        \draw[thick, amber, style=double] (1) -- (2);
        \end{tikzpicture}
\end{subfigure}
\begin{subfigure}{.4\linewidth}
\centering
\begin{tikzpicture}[labelstyle/.style={
       circle, draw=black,
       thick, inner sep=2pt, minimum size=1 mm,
       outer sep=0pt
        }]
       \node (1) at  (0,0) [labelstyle]{$1$};
        \node (2) at (1,0) [labelstyle]{$2$};
        \node (3) at (1,1) [labelstyle]{$3$};
        \node (4) at (2,0) [labelstyle]{$4$};
        \node (5) at (2,1) [labelstyle]{5};
        \node (6) at (3,0) [labelstyle]{$6$};
        \node (7) at (3,1) [labelstyle]{$7$};
        \node (8) at (4,0) [labelstyle]{$8$};
        \node (9) at (4,1) [labelstyle]{$9$};
        \node (10) at (5,0) [labelstyle]{$10$};
        \node (11) at (5,1) [labelstyle]{$11$};

        \draw[thick,red] (2) -- (3);
        \draw[thick,red] (4) -- (5);
        \draw[thick,red] (6) -- (7);
        \draw[thick,red] (8) -- (9);
        \draw[thick,red] (10) -- (11);

        \draw[thick,blue, dashed] (1) -- (3);
        \draw[thick,blue, dashed] (4) -- (6);
        \draw[thick,blue, dashed] (5) -- (7);
        \draw[thick,blue, dashed] (8) -- (10);
        \draw[thick,blue, dashed] (9) -- (11);

        \draw[thick, amber, style=double] (7) -- (9);
        \draw[thick,amber, style=double] (6) -- (11);
        \draw[thick, amber, style=double] (1) -- (2);
        \end{tikzpicture}
\end{subfigure}
\begin{subfigure}{0.4\linewidth}
\begin{tikzpicture}[labelstyle/.style={
       circle, draw=black,
       thick, inner sep=2pt, minimum size=1 mm,
       outer sep=0pt
        }]
       \node (1) at  (0,0) [labelstyle]{$1$};
        \node (2) at (1,0) [labelstyle]{$2$};
        \node (3) at (1,1) [labelstyle]{$3$};
        \node (4) at (2,0) [labelstyle]{$4$};
        \node (5) at (2,1) [labelstyle]{5};
        \node (6) at (3,0) [labelstyle]{$6$};
        \node (7) at (3,1) [labelstyle]{$7$};
        \node (8) at (4,0) [labelstyle]{$8$};
        \node (9) at (4,1) [labelstyle]{$9$};
        \node (10) at (5,0) [labelstyle]{$10$};
        \node (11) at (5,1) [labelstyle]{$11$};

        \draw[thick,red] (2) -- (3);
        \draw[thick,red] (4) -- (5);
        \draw[thick,red] (6) -- (7);
        \draw[thick,red] (8) -- (9);
        \draw[thick,red] (10) -- (11);

        \draw[thick,blue, dashed] (1) -- (3);
        \draw[thick,blue, dashed] (4) -- (6);
        \draw[thick,blue, dashed] (5) -- (7);
        \draw[thick,blue, dashed] (8) -- (10);
        \draw[thick,blue, dashed] (9) -- (11);

        \draw[thick,amber, style=double] (7) -- (9);
        \draw[thick,amber, style=double] (6) to [bend left=30] (10);
        \draw[thick, amber, style=double] (1) -- (2);
        \end{tikzpicture}
\end{subfigure}
\begin{subfigure}{.4\linewidth}
\centering
\begin{tikzpicture}[labelstyle/.style={
       circle, draw=black,
       thick, inner sep=2pt, minimum size=1 mm,
       outer sep=0pt
        }]
       \node (1) at  (0,0) [labelstyle]{$1$};
        \node (2) at (1,0) [labelstyle]{$2$};
        \node (3) at (1,1) [labelstyle]{$3$};
        \node (4) at (2,0) [labelstyle]{$4$};
        \node (5) at (2,1) [labelstyle]{5};
        \node (6) at (3,0) [labelstyle]{$6$};
        \node (7) at (3,1) [labelstyle]{$7$};
        \node (8) at (4,0) [labelstyle]{$8$};
        \node (9) at (4,1) [labelstyle]{$9$};
        \node (10) at (5,0) [labelstyle]{$10$};
        \node (11) at (5,1) [labelstyle]{$11$};

        \draw[thick,red] (2) -- (3);
        \draw[thick,red] (4) -- (5);
        \draw[thick,red] (6) -- (7);
        \draw[thick,red] (8) -- (9);
        \draw[thick,red] (10) -- (11);

        \draw[thick,blue, dashed] (1) -- (3);
        \draw[thick,blue, dashed] (4) -- (6);
        \draw[thick,blue, dashed] (5) -- (7);
        \draw[thick,blue, dashed] (8) -- (10);
        \draw[thick,blue, dashed] (9) -- (11);

        \draw[thick, amber, style=double] (7) -- (9);
        \draw[thick,amber, style=double] (6) -- (8);
        \draw[thick,amber, style=double] (1) -- (2);
        \end{tikzpicture}
\end{subfigure}
\end{center}
\caption{Possibilities for the edge of color $3$ covering vertex $6$, if $c(7,9) = 3$}\label{73_cases}
\end{figure}

Suppose $c(\{7,4\}) = 3$. Up to symmetry, either $c(\{6,5\}) = 3$ or $c(\{6,8\}) = 3$.  If $c(\{6,8\}) = 3$, then $(5,7,4,6,8)$ is a bichromatic $P_4$, so the only possibility is $c(\{6,5\}) = 3$.

\begin{figure}[h]
\begin{center}
\begin{subfigure}{0.4\linewidth}
\centering
\begin{tikzpicture}[labelstyle/.style={
       circle, draw=black,
       thick, inner sep=2pt, minimum size=1 mm,
       outer sep=0pt
        }]
       \node (1) at  (0,0) [labelstyle]{$1$};
        \node (2) at (1,0) [labelstyle]{$2$};
        \node (3) at (1,1) [labelstyle]{$3$};
        \node (4) at (2,0) [labelstyle]{$4$};
        \node (5) at (2,1) [labelstyle]{5};
        \node (6) at (3,0) [labelstyle]{$6$};
        \node (7) at (3,1) [labelstyle]{$7$};
        \node (8) at (4,0) [labelstyle]{$8$};
        \node (9) at (4,1) [labelstyle]{$9$};
        \node (10) at (5,0) [labelstyle]{$10$};
        \node (11) at (5,1) [labelstyle]{$11$};

        \draw[thick,red] (2) -- (3);
        \draw[thick,red] (4) -- (5);
        \draw[thick,red] (6) -- (7);
        \draw[thick,red] (8) -- (9);
        \draw[thick,red] (10) -- (11);

        \draw[thick,blue, dashed] (1) -- (3);
        \draw[thick,blue, dashed] (4) -- (6);
        \draw[thick,blue, dashed] (5) -- (7);
        \draw[thick,blue, dashed] (8) -- (10);
        \draw[thick,blue, dashed] (9) -- (11);

        \draw[thick, amber, style=double] (7) -- (4);
        \draw[thick,amber, style=double] (6) -- (8);
        \draw[thick,amber, style=double] (1) -- (2);
        \end{tikzpicture}
\end{subfigure}
\begin{subfigure}{0.4\linewidth}
\centering
\begin{tikzpicture}[labelstyle/.style={
       circle, draw=black,
       thick, inner sep=2pt, minimum size=1 mm,
       outer sep=0pt
        }]
       \node (1) at  (0,0) [labelstyle]{$1$};
        \node (2) at (1,0) [labelstyle]{$2$};
        \node (3) at (1,1) [labelstyle]{$3$};
        \node (4) at (2,0) [labelstyle]{$4$};
        \node (5) at (2,1) [labelstyle]{5};
        \node (6) at (3,0) [labelstyle]{$6$};
        \node (7) at (3,1) [labelstyle]{$7$};
        \node (8) at (4,0) [labelstyle]{$8$};
        \node (9) at (4,1) [labelstyle]{$9$};
        \node (10) at (5,0) [labelstyle]{$10$};
        \node (11) at (5,1) [labelstyle]{$11$};

        \draw[thick,red] (2) -- (3);
        \draw[thick,red] (4) -- (5);
        \draw[thick,red] (6) -- (7);
        \draw[thick,red] (8) -- (9);
        \draw[thick,red] (10) -- (11);

        \draw[thick,blue, dashed] (1) -- (3);
        \draw[thick,blue, dashed] (4) -- (6);
        \draw[thick,blue, dashed] (5) -- (7);
        \draw[thick,blue, dashed] (8) -- (10);
        \draw[thick,blue, dashed] (9) -- (11);

        \draw[thick, amber, style=double] (7) -- (4);
        \draw[thick,amber, style=double] (6) -- (5);
        \draw[thick,amber, style=double] (1) -- (2);
        \end{tikzpicture}
\end{subfigure}
\end{center}
\caption{Possibilities for the edge of color 3 covering vertex 6, if $c(\{7,4\}) = 3$}
\end{figure}
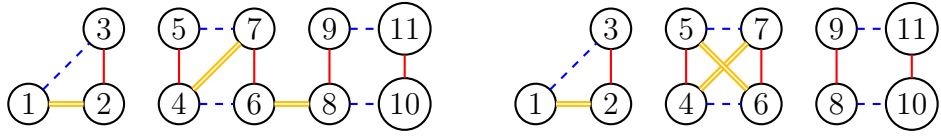

At this point, we have shown that the matching $M_3$ must fall into one of the two cases at the top of  Figure~\ref{reducedCase2}. 
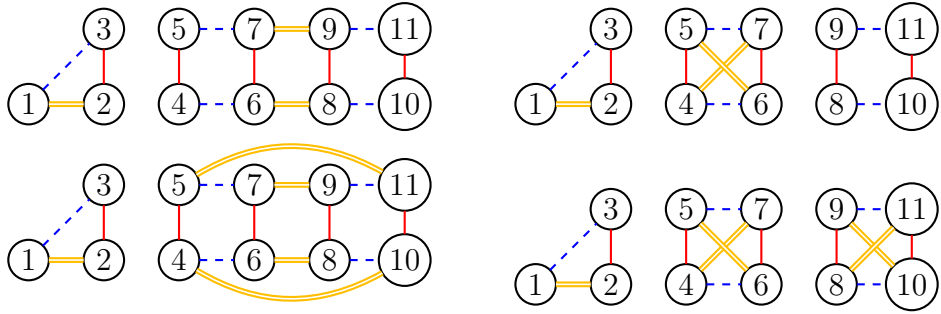
\begin{figure}[h!]
\begin{center}
\begin{subfigure}{.4\linewidth}
\centering
\begin{tikzpicture}[labelstyle/.style={
       circle, draw=black,
       thick, inner sep=2pt, minimum size=1 mm,
       outer sep=0pt
        }]
       \node (1) at  (0,0) [labelstyle]{$1$};
        \node (2) at (1,0) [labelstyle]{$2$};
        \node (3) at (1,1) [labelstyle]{$3$};
        \node (4) at (2,0) [labelstyle]{$4$};
        \node (5) at (2,1) [labelstyle]{5};
        \node (6) at (3,0) [labelstyle]{$6$};
        \node (7) at (3,1) [labelstyle]{$7$};
        \node (8) at (4,0) [labelstyle]{$8$};
        \node (9) at (4,1) [labelstyle]{$9$};
        \node (10) at (5,0) [labelstyle]{$10$};
        \node (11) at (5,1) [labelstyle]{$11$};

        \draw[thick,red] (2) -- (3);
        \draw[thick,red] (4) -- (5);
        \draw[thick,red] (6) -- (7);
        \draw[thick,red] (8) -- (9);
        \draw[thick,red] (10) -- (11);

        \draw[thick,blue, dashed] (1) -- (3);
        \draw[thick,blue, dashed] (4) -- (6);
        \draw[thick,blue, dashed] (5) -- (7);
        \draw[thick,blue, dashed] (8) -- (10);
        \draw[thick,blue, dashed] (9) -- (11);

        \draw[thick, amber, style=double] (7) -- (9);
        \draw[thick,amber, style=double] (6) -- (8);
        \draw[thick,amber, style=double] (1) -- (2);
        \end{tikzpicture}
\end{subfigure}
\begin{subfigure}{0.4\linewidth}
\centering
\begin{tikzpicture}[labelstyle/.style={
       circle, draw=black,
       thick, inner sep=2pt, minimum size=1 mm,
       outer sep=0pt
        }]
       \node (1) at  (0,0) [labelstyle]{$1$};
        \node (2) at (1,0) [labelstyle]{$2$};
        \node (3) at (1,1) [labelstyle]{$3$};
        \node (4) at (2,0) [labelstyle]{$4$};
        \node (5) at (2,1) [labelstyle]{5};
        \node (6) at (3,0) [labelstyle]{$6$};
        \node (7) at (3,1) [labelstyle]{$7$};
        \node (8) at (4,0) [labelstyle]{$8$};
        \node (9) at (4,1) [labelstyle]{$9$};
        \node (10) at (5,0) [labelstyle]{$10$};
        \node (11) at (5,1) [labelstyle]{$11$};

        \draw[thick,red] (2) -- (3);
        \draw[thick,red] (4) -- (5);
        \draw[thick,red] (6) -- (7);
        \draw[thick,red] (8) -- (9);
        \draw[thick,red] (10) -- (11);

        \draw[thick,blue, dashed] (1) -- (3);
        \draw[thick,blue, dashed] (4) -- (6);
        \draw[thick,blue, dashed] (5) -- (7);
        \draw[thick,blue, dashed] (8) -- (10);
        \draw[thick,blue, dashed] (9) -- (11);

        \draw[thick, amber, style=double] (7) -- (4);
        \draw[thick,amber, style=double] (6) -- (5);
        \draw[thick,amber, style=double] (1) -- (2);
        \end{tikzpicture}
\end{subfigure}
\begin{subfigure}{.4\linewidth}
\centering
\begin{tikzpicture}[labelstyle/.style={
       circle, draw=black,
       thick, inner sep=2pt, minimum size=1 mm,
       outer sep=0pt
        }]
       \node (1) at  (0,0) [labelstyle]{$1$};
        \node (2) at (1,0) [labelstyle]{$2$};
        \node (3) at (1,1) [labelstyle]{$3$};
        \node (4) at (2,0) [labelstyle]{$4$};
        \node (5) at (2,1) [labelstyle]{5};
        \node (6) at (3,0) [labelstyle]{$6$};
        \node (7) at (3,1) [labelstyle]{$7$};
        \node (8) at (4,0) [labelstyle]{$8$};
        \node (9) at (4,1) [labelstyle]{$9$};
        \node (10) at (5,0) [labelstyle]{$10$};
        \node (11) at (5,1) [labelstyle]{$11$};

        \draw[thick,red] (2) -- (3);
        \draw[thick,red] (4) -- (5);
        \draw[thick,red] (6) -- (7);
        \draw[thick,red] (8) -- (9);
        \draw[thick,red] (10) -- (11);

        \draw[thick,blue, dashed] (1) -- (3);
        \draw[thick,blue, dashed] (4) -- (6);
        \draw[thick,blue, dashed] (5) -- (7);
        \draw[thick,blue, dashed] (8) -- (10);
        \draw[thick,blue, dashed] (9) -- (11);

        \draw[thick, amber, style=double] (7) -- (9);
        \draw[thick,amber, style=double] (6) -- (8);
        \draw[thick,amber, style=double] (1) -- (2);
        \draw[thick, amber, style=double] (5) to [bend left=30] (11);
        \draw[thick,amber, style=double] (4) to [bend right=30] (10);
        \end{tikzpicture}
\end{subfigure}
\begin{subfigure}{0.4\linewidth}
\centering
\begin{tikzpicture}[labelstyle/.style={
       circle, draw=black,
       thick, inner sep=2pt, minimum size=1 mm,
       outer sep=0pt
        }]
       \node (1) at  (0,0) [labelstyle]{$1$};
        \node (2) at (1,0) [labelstyle]{$2$};
        \node (3) at (1,1) [labelstyle]{$3$};
        \node (4) at (2,0) [labelstyle]{$4$};
        \node (5) at (2,1) [labelstyle]{5};
        \node (6) at (3,0) [labelstyle]{$6$};
        \node (7) at (3,1) [labelstyle]{$7$};
        \node (8) at (4,0) [labelstyle]{$8$};
        \node (9) at (4,1) [labelstyle]{$9$};
        \node (10) at (5,0) [labelstyle]{$10$};
        \node (11) at (5,1) [labelstyle]{$11$};

        \draw[thick,red] (2) -- (3);
        \draw[thick,red] (4) -- (5);
        \draw[thick,red] (6) -- (7);
        \draw[thick,red] (8) -- (9);
        \draw[thick,red] (10) -- (11);

        \draw[thick,blue, dashed] (1) -- (3);
        \draw[thick,blue, dashed] (4) -- (6);
        \draw[thick,blue, dashed] (5) -- (7);
        \draw[thick,blue, dashed] (8) -- (10);
        \draw[thick,blue, dashed] (9) -- (11);

        \draw[thick, amber, style=double] (7) -- (4);
        \draw[thick,amber, style=double] (6) -- (5);
        \draw[thick,amber, style=double] (1) -- (2);
        \draw[thick, amber, style=double] (10) -- (9);
        \draw[thick,amber, style=double] (11) -- (8);
        
        \end{tikzpicture}
\end{subfigure}
\end{center}
\caption{Subcases of Case 2}\label{reducedCase2}
\end{figure}
Note that we can determine the remaining edges in $M_3$ in each case, since we cannot form bichromatic $P_4$'s. 
The complete pictures for $M_1,M_2$ and $M_3$ are given at the bottom of 
Figure~\ref{reducedCase2}.

In the first subcase, the four-cycle property implies that $M_4$ has a perfect matching of $\{1,2,3,5,6,10\}$ as a submatching. Since $4$ is not covered by $M_4$, the remaining edges of $M_4$ form a perfect matching of $\{7,8,9,11\}$. Since $9$ is already adjacent to each of $7,8,11$, this is impossible. In the second subcase, the four-cycle property implies that $M_4$ has a perfect matching of $\{1,2,3,5,6,7\}$ as a submatching. Since $4$ is not covered by $M_4$, the remaining edges in $M_4$ must form a perfect matching of $\{8,9,10,11\}$. However, all edges with two endpoints in this set have already been colored by $M_1,M_2$ and $M_3$. So $M_4$ cannot form a perfect matching of these vertices and the second subcase is not possible. All possible cases have been eliminated, and so the coloring $c$ does not exist.
\end{proof}
\section{Concluding Remarks} \label{conclusion}

We demonstrate that the behavior of $\mathrm{ex^*}(n,B_{t,3})$ depends not only upon the parity of $t$, but (at least) upon whether $t + 2$ is a power of $2$. In some cases when $t+2$ is not a power of $2$, we are able to show that $\mathrm{ex^*}(n,B_{t,3}) = \frac{t}{2}n + O(1)$. It is tempting to conjecture that  $\mathrm{ex^*}(n,B_{t,3}) = \frac{t}{2}n + O(1)$ whenever $t+2$ is not a power of $2$, or at least that $\mathrm{ex^*}(n,B_{t,3}) = \frac{t}{2}n + O(1)$ whenever $t \not\equiv 2 \pmod{4}$. The main barrier to proving (or disproving) such a conjecture seems to be in finding extremal constructions; the algebraic colorings used to generate the constructions in Sections~\ref{power of 2 construction} and \ref{power of 3 construction} will necessarily only apply for a sparse set of $t$ values. 

We also remark that stability will not hold in general for these problems; that is, we cannot expect that there is a unique extremal construction achieving $\mathrm{ex^*}(n,B_{t,3})$. Indeed, since $2^3 = 3^2 - 1$, the (non-isomorphic) constructions presented in Sections~\ref{power of 2 construction} and \ref{power of 3 construction} furnish two extremal graphs for $\mathrm{ex^*}(n,B_{8,3})$. On the other hand, Lemma~\ref{extremal space restriction} and Proposition~\ref{powers of 2} imply that when $t + 2 = 2^s$, then the unique extremal graph for $\mathrm{ex^*}(n,B_{t,3})$ is a disjoint union of $K_{t+2}$-copies, each receiving an optimal edge-coloring. It would be interesting to characterize other instances where a unique extremal construction does (or does not) exist for $\mathrm{ex^*}(n,B_{t,3})$; such investigation might also shed light on the problem of finding good lower-bound constructions for $\mathrm{ex^*}(n,B_{t,3})$ which generalize more easily than the ones provided here.

In light of our results, the first case in which $\mathrm{ex^*}(n,B_{t,3})$ is unresolved is $B_{10,3}$. In Proposition~\ref{B10structure}, we demonstrate that $K_{11}$ cannot be colored to avoid a rainbow $B_{10,3}$ copy, so either $\mathrm{ex}(n,B_{10,3}) = \frac{9}{2}n + O(1)$ (which is achievable by taking disjoint copies of $K_{10}$ with an arbitrary proper edge-coloring) or else $\mathrm{ex^*}(n,B_{10,3})$ is not achieved by cliques. Either outcome would be interesting. In the first case, we would have $\mathrm{ex^*}(n,B_{10,3}) = \mathrm{ex^*}(n,B_{9,3}) + O(1)$. In the second, it would follow that there is no ``unified'' extremal construction for $B_{t,3}$: cliques furnish the unique extremal construction when $t = 2^s - 2$, but would be provably suboptimal for $B_{10,3}$.

\section{Acknowledgments}

This work was started at the 2024 Graduate Research Workshop in Combinatorics. The authors would like to thank the workshop hosts, organizers, and faculty mentors. We would in particular like to thank Puck Rombach and Shira Zerbib for useful conversations and support during the early stages of this project. 

\bibliographystyle{abbrvurl}
\bibliography{master.bib}

\begin{thebibliography}{1}

\bibitem{bednar2022rainbow}
V.~Bednar and N.~Bushaw.
\newblock Rainbow {T}ur{\'a}n {M}ethods for {T}rees.
\newblock {\em arXiv preprint arXiv:2203.13765}, 2022.

\bibitem{ErdosGallai}
P.~Erd\H{o}s and T.~Gallai.
\newblock On maximal paths and circuits of graphs.
\newblock {\em Act. Math. Hungar.}, 10:337--356, 1959.

\bibitem{ErSi}
P.~Erd\H{o}s and M.~Simonovits.
\newblock A limit theorem in graph theory.
\newblock {\em Studia Sci. Math. Hungar.}, 1:51--57, 1966.

\bibitem{ErSt}
P.~Erd{\H o}s and A.~H. Stone.
\newblock On the structure of linear graphs.
\newblock {\em Bull. Amer. Math. Soc.}, 52:1087--1091, 1946.
\newblock \href {https://doi.org/10.1090/S0002-9904-1946-08715-7} {\path{doi:10.1090/S0002-9904-1946-08715-7}}.

\bibitem{halfpaprainbowp5}
A.~Halfpap.
\newblock The rainbow {T}ur\'an number of {$P_5$}.
\newblock {\em Australas. J. Combin.}, 87 (3):403--422, 2023.

\bibitem{JoPaSa}
D.~Johnston, C.~Palmer, and A.~Sarkar.
\newblock Rainbow {T}ur\'{a}n problems for paths and forests of stars.
\newblock {\em Electron. J. Combin.}, 24(1):Paper No. 1.34, 15, 2017.
\newblock \href {https://doi.org/10.37236/6430} {\path{doi:10.37236/6430}}.

\bibitem{JoRo}
D.~Johnston and P.~Rombach.
\newblock Lower bounds for rainbow {T}ur\'{a}n numbers of paths and other trees.
\newblock {\em Australas. J. Combin.}, 78:61--72, 2020.

\bibitem{KMSV}
P.~Keevash, D.~Mubayi, B.~Sudakov, and J.~Verstra\"{e}te.
\newblock Rainbow {T}ur\'{a}n problems.
\newblock {\em Combin. Probab. Comput.}, 16(1):109--126, 2007.
\newblock \href {https://doi.org/10.1017/S0963548306007760} {\path{doi:10.1017/S0963548306007760}}.

\bibitem{Kirkman}
T.~Kirkman.
\newblock On a problem in combinations.
\newblock {\em Cambridge and Dublin Math. J.}, 2:191--204, 1847.

\end{thebibliography}

\end{document}